\documentclass[11pt]{article}
\renewcommand*{\today}{\number\day\space\ifcase\month\or January\or
  February\or March\or April\or May\or June\or July\or August\or
  September\or October\or November\or December\fi\space\number\year}

\usepackage{amssymb}
\usepackage{amsmath}
\usepackage{amsthm}
\usepackage{a4wide}
\usepackage{color}
\usepackage{url}

\usepackage{graphics}
\usepackage[pdftex]{graphicx}

\newcommand{\arrc}{\overrightarrow{C}}

\newtheorem{thm}{Theorem}
\newtheorem{lem}[thm]{Lemma}
\newtheorem{prop}[thm]{Proposition}
\newtheorem*{ques}{Question}

\begin{document}

\title{\textbf{Connectedness of Strong $k$-Colour Graphs}}
\author{Somkiat Trakultraipruk\\[2mm]
  \emph{Department of Mathematics}\\
  \emph{London School of Economics and Political Science, London, U.K.}}
\maketitle

\begin{abstract}
  \noindent
  For a positive integer~$k$ and a graph~$G$, we consider proper
  vertex-colourings of~$G$ with~$k$ colours in which all~$k$ colours are
  actually used. We call such a colouring a strong \mbox{$k$-colouring}.
  The strong $k$-colour graph of~$G$, $S_k(G)$, is the graph that has all
  the strong $k$-colourings of~$G$ as its vertex set, and two colourings
  are adjacent in $S_k(G)$ if they differ in colour on only one vertex
  of~$G$. In this paper, we show some results related to the question\,:
  For what~$G$ and~$k$ is $S_k(G)$ connected\,?

\end{abstract}

\section{Introduction}\label{sec1}

  Throughout this paper a graph~$G$ is finite, simple, and loopless, and we
  also usually assume that~$G$ is connected. Most of our terminology and
  notation will be standard and can be found in any textbook on graph theory
  such as~\cite{RDiestel} and~\cite{DBWest}. For a positive integer~$k$ and a
  graph~$G$, the \emph{$k$-colour graph} of~$G$, denoted $C_k(G)$, is the
  graph that has the proper $k$-vertex-colourings of~$G$ as its vertex set,
  and two such colourings are joined by an edge in~$C_k(G)$ if they differ in
  colour on only one vertex of~$G$.

  We now introduce a subgraph of $C_k(G)$, called the \emph{strong $k$-colour
  graph} of~$G$, denoted $S_k(G)$. Its vertex set contains only proper
  $k$-colourings in which all~$k$ colours actually appear, and we call such a
  colouring a \emph{strong $k$-colouring}.

  Questions regarding the connectivity of a $k$-colour graph have
  applications in reassignment problems of the channels used in cellular
  networks; see, e.g.,~\cite{freq-ass,JHan,RJac}. For some applications, it
  is required that all channels in a range are actually used. Such a
  labelling is sometimes called a ``no-hole'' or ``consecutive'' labelling;
  see, e.g.,~\cite{no-hole,cons}. In terms of colourings, this corresponds to
  a strong $k$-colouring. And asking questions about the possibility to
  reassign channels in a cellular network can be done in such a way that all
  available channels are actually used. These problems can be expressed in
  finding paths in the strong $k$-colour graph.

\begin{figure}[ht]
  \centering
  \includegraphics[height=0.24\textheight]{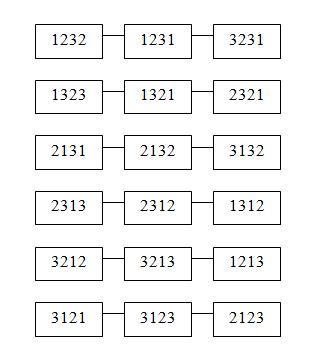}
  \caption{The strong colour graph $S_3(P_4)$.}
  \label{f1}
\end{figure}

\begin{figure}[ht]
  \centering
  \includegraphics[height=0.34\textheight]{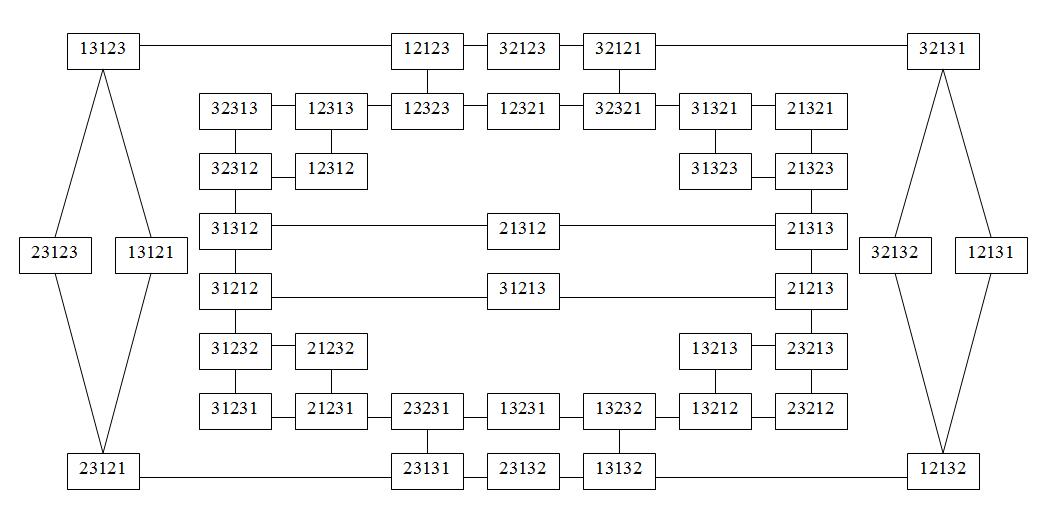}
  \caption{The strong colour graph $S_3(P_5)$.}
  \label{f2}
\end{figure}

  Questions related to the connectivity of a $k$-colour graph have been
  studied extensively\,: \cite{vertex-colouring,mixing,YMCAL,OMPH}. In this
  note we initiate similar research on the connectivity of strong $k$-colour
  graphs\,: Given a positive integer~$k$ and a graph~$G$, is $S_k(G)$
  connected? As an example, in Figures~\ref{f1} and~\ref{f2}, we show the
  strong 3-colour graph of the paths with~4 and~5 vertices, respectively. One
  is not connected while the other is connected.

\medskip
  We usually use lower case Greek letters $\alpha,\beta,\gamma,\ldots$ to
  denote specific colourings, and lower case Latin $a,b,c,\ldots$ to denote
  specific colours.

  To avoid trivial cases, we will always assume that~$k$ is greater than or
  equal to the chromatic number of~$G$, and the number of vertices of~$G$ is
  at least $k+1$.

\section{General Results}\label{sec2}

  For all $k\ge2$ and $m,n\ge1$, let~$\alpha$ be a strong
  $k$-vertex-colouring of a complete bipartite graph $K_{m,n}$.
  Each colour that appears in one part of the partition cannot be used
  in the other part. Now we choose one colour from each part and recolour the
  graph by swapping these two colours on each vertex coloured with one of
  these colours. Let~$\beta$ be the resulting colouring, so~$\beta$ is strong
  as well. It is easy to see that there is no path in $S_k(K_{m,n})$
  from~$\alpha$ to~$\beta$. Thus $S_k(K_{m,n})$ is not connected for all
  $k\ge2$ and $m,n\ge1$. This gives the following lemma.

\begin{lem}\label{Z.8}
  Let~$G$ be a connected, $k$-colourable graph such that $S_k(G)$ is
  connected, with $k\ge2$. Then $|V(G)|\ge k+1$, and~$G$ does not contain a
  complete bipartite graph as a spanning subgraph.
\end{lem}

\begin{thm}\label{Z.6}
  Let~$G$ be a connected, $k$-colourable graph such that $S_k(G)$ is
  connected, with $k\ge2$. Suppose the graph~$G^*$ is obtained from~$G$ by
  adding a new vertex~$v^*$ and joining it to~$j$ vertices in $V(G)$, with
  $1\le j\le k-2$. Then $S_k(G^{*})$ is connected.
\end{thm}

\noindent
  We will show in the next section that the strong $3$-colour graph of
  the $n$-vertex path, $S_3(P_n)$, is connected
  if and only if $n\ge5$. Now add
  a new vertex and join it to the first and the last vertex of
  the path, forming to an $(n+1)$-vertex cycle.
  We will show in Section~\ref{sec4} that the strong $3$-colour graph
  of the $n$-vertex cycle, $S_3(C_n)$, is not connected for all~$n$.
  This example shows that the restriction $j\le k-2$ of Theorem~\ref{Z.6} is
  optimal.
  
\begin{proof}[\textbf{Proof of Theorem~\ref{Z.6}}]
  Let~$\alpha^*$ and~$\beta^*$ be strong $k$-colourings of~$G^*$. We show
  that there always exists a walk in $S_k(G^*)$ from~$\alpha^*$
  to~$\beta^*$. We say that a colouring of~$G^*$ is \emph{good} if all~$k$
  colours appear on~$V(G)$.

  First suppose that~$\alpha^*$ and~$\beta^*$ are good. By ignoring the
  vertex~$v^*$, let~$\alpha$ and~$\beta$ be the strong $k$-colourings
  of~$G$ obtained from~$\alpha^*$ and~$\beta^*$, respectively. Since
  $S_k(G)$ is connected, there is a path from~$\alpha$ to~$\beta$ in
  $S_k(G)$. We just follow the recolouring steps of that path to form a
  walk from~$\alpha^*$ to~$\beta^*$ in $S_k(G^*)$. The only extra steps
  happen when we want to recolour a neighbour~$u$ of~$v^*$ to the same
  colour as~$v^*$. Since $d_{G^*}(v^*)=j\le k-2$, we can always
  recolour~$v^*$ to
  a colour different from any of the colours appearing in its
  neighbourhood and its current colour. After recolouring~$v^*$, we can
  recolour~$u$, and continue the walk. This walk in $S_k(G^*)$ finishes in
  a colouring in which the vertices in~$V(G)$ have the same colour as they
  have in~$\beta^*$. If necessary, we can do one recolouring of~$v^*$ to
  its colour in~$\beta^*$, completing the walk in $S_k(G^*)$
  from~$\alpha^*$ to~$\beta^*$.

  If~$\alpha^*$ is not good, then below we show that we can always find a
  path in $S_k(G^*)$ from~$\alpha^*$ to some good colouring (\,and if
  necessary, we do the same for~$\beta^*$\,). Together with the method
  described in the previous paragraph, this completes the proof.

  So we now assume that in~$\alpha^*$ every vertex of~$G$ has received one
  of $k-1$ colours while~$v^*$ has the remaining colour. Let~$W$ be the set
  of vertices in~$V(G)$ that do not have a unique colour in~$G$ for the
  colouring~$\alpha^*$. Since $|V(G)|\ge k+1$, $W$ is not empty.

  \medskip\noindent
  \textbf{Case 1}:\quad There is a vertex $w\in W$ not adjacent
  to~$v^*$.\\*[\parskip]
  By the definition of~$W$, there is a vertex $w'\in W$ such that~$w$
  and~$w'$ have the same colour in~$\alpha^*$. Then we recolour~$w$ to the
  same colour as~$v^*$. The resulting colouring is good.

  \medskip\noindent
  \textbf{Case 2}:\quad All vertices in~$W$ are adjacent
  to~$v^*$.\\*[\parskip]
  Additionally, define $U=N(v^*)\setminus W$. and $X=V(G)\setminus N(v^*)$.
  Note that all vertices in~$X$ have a unique colour in~$\alpha^*$.

  \medskip\noindent \textbf{Subcase 2.1}:\quad There is a vertex $x\in X$
  that is not adjacent to some vertex $w\in W$.\\*[\parskip]
  Again, there is a vertex $w'\in W$ such that~$w$ and~$w'$ have the same
  colour in~$\alpha^*$. Then we first recolour~$w$ to the same colour
  as~$x$, and then recolour~$x$ to the same colour as~$v^*$. Again, this
  gives a good colouring.

  \medskip\noindent

  \textbf{Subcase 2.2}:\quad Every vertex in~$X$ is adjacent to every
  vertex in~$W$.\\*[\parskip]
  Because of Lemma~\ref{Z.8}, $U$ is not empty (\,otherwise, the pair
  $(X,W)$ would form the parts of a spanning complete bipartite subgraph
  of~$G$\,). Suppose there is some vertex $u\in U$ that is not adjacent to
  some vertex $w\in W$ and not adjacent to some vertex $x\in X$. Then we
  can recolour~$w$ to the colour of~$u$ (\,this is possible since there is
  another vertex $w'\in W$ with the same colour as~$w$\,). Then
  recolour~$u$ to the same colour as~$x$, and lastly~$x$ to the same colour
  as~$v^*$. It is easy to check that the remaining colouring is good.

  \medskip
  So we are left with the case that each vertex in~$U$ is adjacent to every
  vertex in~$W$ or to every vertex in~$X$. Let~$U_W$ be the set of vertices
  in~$U$ that are adjacent to every vertex in~$W$, and
  $U_X=U\setminus U_W$. Then the pair $(X\cup U_W,W\cup U_X)$ forms the
  parts of a spanning complete bipartite subgraph of~$G$. Because of
  Lemma~\ref{Z.8}, this contradicts that $S_k(G)$ is connected.
\end{proof}

\begin{thm}\label{Z.5}
  Let~$G$ be a connected $k$-colourable graph so that $S_k(G)$ is
  connected, with $k\ge2$. Let~$v$ be a vertex of~$G$ with neighbourhood
  $N(v)$. Suppose the graph~$G^*$ is obtained from~$G$ by adding a new
  vertex~$v^*$ and joining~$v^*$ to the vertices in~$N^*$ for some
  $N^*\subseteq N(v)$, $N^*\ne\varnothing$. Then $S_k(G^{*})$ is connected.
\end{thm}

\begin{proof}
  Let~$\alpha^*$ and~$\beta^*$ be strong $k$-colourings of~$G^*$. We show
  that there always exists a walk in $S_k(G^*)$ from~$\alpha^*$
  to~$\beta^*$. We say that a colouring of~$G^*$ is \emph{good} if~$v$
  and~$v^*$ are labelled with the same colour.

  First suppose that~$\alpha^*$ and~$\beta^*$ are good. By ignoring the
  vertex~$v^*$, let~$\alpha$ and~$\beta$ be the strong $k$-colourings
  of~$G$ obtained from~$\alpha^*$ and~$\beta^*$, respectively. Since
  $S_k(G)$ is connected, there is a path from~$\alpha$ to~$\beta$ in
  $S_k(G)$. We just follow the recolouring steps of that path to form a
  walk from~$\alpha^*$ to~$\beta^*$ in $S_k(G^*)$. The only extra step
  happens when we recolour~$v$. In the next step we immediately
  recolour~$v^*$ to the same colour as~$v$ just received. It is easy to
  check that all these recolourings are allowed and give a walk in
  $S_k(G^*)$ from~$\alpha^*$ to~$\beta^*$, completing the proof.

  Assume that~$\alpha^*$ is not good. Below we show that we always can find
  a path in $S_k(G^*)$ from~$\alpha^*$ to some good colouring (\,and if
  necessary, we do the same for~$\beta^*$\,). Together with the method
  described in the previous paragraph, this completes the proof.

  If there is a vertex $u\in V(G)\setminus\{v\}$ with the same colour
  as~$v^*$, then we just recolour~$v^*$ to the same colour as~$v$. This
  gives a good colouring.

  So we now suppose that in~$\alpha^*$ every vertex of~$G$ has received one
  of $k-1$ colours while~$v^*$ has the remaining colour. Remind that~$v$
  and~$v^*$ received different colours in~$\alpha^*$. Let~$W$ be the set of
  vertices in~$V(G)$ that did not receive a unique colour in~$G$ for the
  colouring~$\alpha^*$. Since $|V(G)|\ge k+1$, $W$ is not empty.

  \medskip\noindent
  \textbf{Case 1}:\quad There is a vertex $w\in W$ not adjacent
  to~$v^*$.\\*[\parskip]
  By the definition of~$W$, there is a vertex $w'\in W$ such that~$w$
  and~$w'$ have the same colour in~$\alpha^*$. Hence we can recolour~$w$ to
  the same colour as~$v^*$, and then recolour~$v^*$ to the same colour
  as~$v$. The resulting colouring is good.

  \medskip\noindent
  \textbf{Case 2}:\quad All vertices in~$W$ are adjacent
  to~$v^*$.\\*[\parskip]
  Additionally, define $U=N(v^*)\setminus W$. and $X=V(G)\setminus N(v^*)$.
  Note that all vertices in~$X$ have a unique colour in~$\alpha^*$.

  \medskip\noindent
  \textbf{Subcase 2.1}:\quad There is a vertex $x\in X$ that is not
  adjacent to some vertex $w\in W$.\\*[\parskip]
  Again, there is a vertex $w'\in W$ such that~$w$ and~$w'$ have the same
  colour in~$\alpha^*$. Then we first recolour~$w$ to the same colour
  as~$x$. Then recolour~$x$ to the same colour as~$v^*$, and lastly
  recolour~$v^*$ to the same colour as~$v$. Again, this gives a good
  colouring.

  \medskip\noindent
  \textbf{Subcase 2.2}:\quad Every vertex in~$X$ is adjacent to every
  vertex in~$W$.\\*[\parskip]
  Because of Lemma~\ref{Z.8}, $U$ is not empty (\,otherwise, the pair
  $(X,W)$ would form the parts of a spanning complete bipartite subgraph
  of~$G$\,). Suppose there is some vertex $u\in U$ that is not adjacent to
  some vertex $w\in W$ and not adjacent to some vertex $x\in X$. Then we
  can recolour~$w$ to the colour of~$u$ (\,this is possible since there is
  another vertex $w'\in W$ with the same colour as~$w$\,). Then
  recolour~$u$ to the same colour as~$x$, ~$x$ to the same colour as~$v^*$,
  and lastly recolour~$v^*$ to the same colour as~$v$. It is easy to check
  that the remaining colouring is good.

  \medskip
  So we are left with the case that each vertex in~$U$ is adjacent to every
  vertex in~$W$ or to every vertex in~$X$. Let~$U_W$ be the set of vertices
  in~$U$ that are adjacent to every vertex in~$W$, and
  $U_X=U\setminus U_W$. Then the pair $(X\cup U_W,W\cup U_X)$ forms the
  parts of a spanning complete bipartite subgraph of~$G$. Because of
  Lemma~\ref{Z.8}, this contradicts that $S_k(G)$ is connected.
\end{proof}

\noindent
  It is easy to see that in the normal colour graph~$C_k(G)$ there always is
  a path from any proper $k$-vertex-colouring to some strong
  $k$-vertex-colouring. This shows the following.

\begin{lem} \label{Z.10}
  If $S_k(G)$ is connected, then $C_k(G)$ is also connected.
\end{lem}

\section{The Strong $k$-Colour Graph of Paths}\label{sec3}

  In this section, we prove that the strong $k$-colour graph of a path
  with~$n$ vertices, $S_k(P_n)$, is connected if and only if $k\ge3$,
  $n\ge5$, and $n\ge k+1$.

  First, suppose we colour a path~$P_n$, $n\ge2$, with two colours. It is
  easy to see that there are only two strong 2-vertex-colourings
  of~$P_n$, and they are not adjacent in $S_2(P_n)$. Thus $S_2(P_n)$ is not
  connected for all $n\ge2$.

  For $k=3$, we have already seen in Figures~\ref{f1} and~\ref{f2} that
  $S_3(P_4)$ is not connected, but $S_3(P_5)$ is connected.
  It is somewhat more work to show that $S_4(P_5)$ is connected.

\begin{prop}\label{B.6}
  The strong colour graph $S_4(P_5)$ is connected.
\end{prop}

\begin{proof}
  In any strong 4-colouring of~$P_5$, there are only two vertices with the
  same colour. Let~$\alpha$ be a strong 4-colouring of $P_5=v_1v_2...v_5$.
  We call~$\alpha$ an \emph{$a$-standard colouring} if
  $\alpha(v_1)=\alpha(v_5)=a$.

  \begin{figure}[ht]
    \centering
    \unitlength0.25mm
    \begin{picture}(130,40)(0,5)
      \put(25,20){\circle*{5}}
      \put(45,20){\circle*{5}}
      \put(65,20){\circle*{5}}
      \put(85,20){\circle*{5}}
      \put(105,20){\circle*{5}}
      \put(25,20){\line(1,0){80}}
      \put(25,33){\makebox(0,0)[b]{$a$}}
      \put(45,33){\makebox(0,0)[b]{$b$}}
      \put(65,33){\makebox(0,0)[b]{$c$}}
      \put(85,33){\makebox(0,0)[b]{$d$}}
      \put(105,33){\makebox(0,0)[b]{$a$}}
      \put(25,5){\makebox(0,0)[t]{$v_1$}}
      \put(105,5){\makebox(0,0)[t]{$v_5$}}
    \end{picture}

    \caption{An $a$-standard colouring.}
    \label{f3}
  \end{figure}
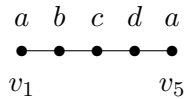

  We will prove the proposition by combining one or more of the following
  three steps.

  \medskip\noindent
  \textbf{Step 1}:\quad There is a path from $(a,b,c,d,a)$ to any other
  $a$-standard colouring.\\*[\parskip]
  We will first show that there is a path from $(a,b,c,d,a)$ to
  $(a,c,b,d,a)$\,:

  \begin{center}
    \unitlength0.20mm
    \begin{picture}(120,40)(5,15)
      \put(25,20){\circle*{5}}
      \put(45,20){\circle*{5}}
      \put(65,20){\circle*{5}}
      \put(85,20){\circle*{5}}
      \put(105,20){\circle*{5}}
      \put(25,20){\line(1,0){80}}
      \put(25,33){\makebox(0,0)[b]{$a$}}
      \put(45,33){\makebox(0,0)[b]{$b$}}
      \put(65,33){\makebox(0,0)[b]{$c$}}
      \put(85,33){\makebox(0,0)[b]{$d$}}
      \put(105,33){\makebox(0,0)[b]{$a$}}
    \end{picture}
    \begin{picture}(25,40)(0,15)
      \put(0,20){\vector(1,0){25}}
    \end{picture}
    \begin{picture}(120,40)(5,15)
      \put(25,20){\circle*{5}}
      \put(45,20){\circle*{5}}
      \put(65,20){\circle*{5}}
      \put(85,20){\circle*{5}}
      \put(105,20){\circle*{5}}
      \put(25,20){\line(1,0){80}}
      \put(25,33){\makebox(0,0)[b]{$a$}}
      \put(45,33){\makebox(0,0)[b]{$b$}}
      \put(65,33){\makebox(0,0)[b]{$c$}}
      \put(85,33){\makebox(0,0)[b]{$d$}}
      \put(105,33){\makebox(0,0)[b]{$b$}}
    \end{picture}
    \begin{picture}(25,40)(0,15)
      \put(0,20){\vector(1,0){25}}
    \end{picture}
    \begin{picture}(120,40)(5,15)
      \put(25,20){\circle*{5}}
      \put(45,20){\circle*{5}}
      \put(65,20){\circle*{5}}
      \put(85,20){\circle*{5}}
      \put(105,20){\circle*{5}}
      \put(25,20){\line(1,0){80}}
      \put(25,33){\makebox(0,0)[b]{$a$}}
      \put(45,33){\makebox(0,0)[b]{$d$}}
      \put(65,33){\makebox(0,0)[b]{$c$}}
      \put(85,33){\makebox(0,0)[b]{$d$}}
      \put(105,33){\makebox(0,0)[b]{$b$}}
    \end{picture}
    \begin{picture}(25,40)(0,15)
      \put(0,20){\vector(1,0){25}}
    \end{picture}
    \begin{picture}(120,40)(5,15)
      \put(25,20){\circle*{5}}
      \put(45,20){\circle*{5}}
      \put(65,20){\circle*{5}}
      \put(85,20){\circle*{5}}
      \put(105,20){\circle*{5}}
      \put(25,20){\line(1,0){80}}
      \put(25,33){\makebox(0,0)[b]{$a$}}
      \put(45,33){\makebox(0,0)[b]{$d$}}
      \put(65,33){\makebox(0,0)[b]{$c$}}
      \put(85,33){\makebox(0,0)[b]{$a$}}
      \put(105,33){\makebox(0,0)[b]{$b$}}
    \end{picture}
    \begin{picture}(25,40)(0,15)
      \put(0,20){\vector(1,0){25}}
    \end{picture}
    \begin{picture}(120,40)(5,15)
      \put(25,20){\circle*{5}}
      \put(45,20){\circle*{5}}
      \put(65,20){\circle*{5}}
      \put(85,20){\circle*{5}}
      \put(105,20){\circle*{5}}
      \put(25,20){\line(1,0){80}}
      \put(25,33){\makebox(0,0)[b]{$c$}}
      \put(45,33){\makebox(0,0)[b]{$d$}}
      \put(65,33){\makebox(0,0)[b]{$c$}}
      \put(85,33){\makebox(0,0)[b]{$a$}}
      \put(105,33){\makebox(0,0)[b]{$b$}}
    \end{picture}
    \begin{picture}(25,40)(0,15)
      \put(0,20){\vector(1,0){25}}
    \end{picture}
    \begin{picture}(120,40)(5,15)
      \put(25,20){\circle*{5}}
      \put(45,20){\circle*{5}}
      \put(65,20){\circle*{5}}
      \put(85,20){\circle*{5}}
      \put(105,20){\circle*{5}}
      \put(25,20){\line(1,0){80}}
      \put(25,33){\makebox(0,0)[b]{$c$}}
      \put(45,33){\makebox(0,0)[b]{$d$}}
      \put(65,33){\makebox(0,0)[b]{$b$}}
      \put(85,33){\makebox(0,0)[b]{$a$}}
      \put(105,33){\makebox(0,0)[b]{$b$}}
    \end{picture}
    \begin{picture}(25,40)(0,15)
      \put(0,20){\vector(1,0){25}}
    \end{picture}
    \begin{picture}(120,40)(5,15)
      \put(25,20){\circle*{5}}
      \put(45,20){\circle*{5}}
      \put(65,20){\circle*{5}}
      \put(85,20){\circle*{5}}
      \put(105,20){\circle*{5}}
      \put(25,20){\line(1,0){80}}
      \put(25,33){\makebox(0,0)[b]{$c$}}
      \put(45,33){\makebox(0,0)[b]{$d$}}
      \put(65,33){\makebox(0,0)[b]{$b$}}
      \put(85,33){\makebox(0,0)[b]{$a$}}
      \put(105,33){\makebox(0,0)[b]{$c$}}
    \end{picture}
    \begin{picture}(25,40)(0,15)
      \put(0,20){\vector(1,0){25}}
    \end{picture}
    \begin{picture}(120,40)(5,15)
      \put(25,20){\circle*{5}}
      \put(45,20){\circle*{5}}
      \put(65,20){\circle*{5}}
      \put(85,20){\circle*{5}}
      \put(105,20){\circle*{5}}
      \put(25,20){\line(1,0){80}}
      \put(25,33){\makebox(0,0)[b]{$a$}}
      \put(45,33){\makebox(0,0)[b]{$d$}}
      \put(65,33){\makebox(0,0)[b]{$b$}}
      \put(85,33){\makebox(0,0)[b]{$a$}}
      \put(105,33){\makebox(0,0)[b]{$c$}}
    \end{picture}
    \begin{picture}(25,40)(0,15)
      \put(0,20){\vector(1,0){25}}
    \end{picture}
    \begin{picture}(120,40)(5,15)
      \put(25,20){\circle*{5}}
      \put(45,20){\circle*{5}}
      \put(65,20){\circle*{5}}
      \put(85,20){\circle*{5}}
      \put(105,20){\circle*{5}}
      \put(25,20){\line(1,0){80}}
      \put(25,33){\makebox(0,0)[b]{$a$}}
      \put(45,33){\makebox(0,0)[b]{$d$}}
      \put(65,33){\makebox(0,0)[b]{$b$}}
      \put(85,33){\makebox(0,0)[b]{$d$}}
      \put(105,33){\makebox(0,0)[b]{$c$}}
    \end{picture}
    \begin{picture}(25,40)(0,15)
      \put(0,20){\vector(1,0){25}}
    \end{picture}
    \begin{picture}(120,40)(5,15)
      \put(25,20){\circle*{5}}
      \put(45,20){\circle*{5}}
      \put(65,20){\circle*{5}}
      \put(85,20){\circle*{5}}
      \put(105,20){\circle*{5}}
      \put(25,20){\line(1,0){80}}
      \put(25,33){\makebox(0,0)[b]{$a$}}
      \put(45,33){\makebox(0,0)[b]{$c$}}
      \put(65,33){\makebox(0,0)[b]{$b$}}
      \put(85,33){\makebox(0,0)[b]{$d$}}
      \put(105,33){\makebox(0,0)[b]{$c$}}
    \end{picture}
    \begin{picture}(25,40)(0,15)
      \put(0,20){\vector(1,0){25}}
    \end{picture}
    \begin{picture}(120,40)(5,15)
      \put(25,20){\circle*{5}}
      \put(45,20){\circle*{5}}
      \put(65,20){\circle*{5}}
      \put(85,20){\circle*{5}}
      \put(105,20){\circle*{5}}
      \put(25,20){\line(1,0){80}}
      \put(25,33){\makebox(0,0)[b]{$a$}}
      \put(45,33){\makebox(0,0)[b]{$c$}}
      \put(65,33){\makebox(0,0)[b]{$b$}}
      \put(85,33){\makebox(0,0)[b]{$d$}}
      \put(105,33){\makebox(0,0)[b]{$a$}}
    \end{picture}
  \end{center}

  By symmetry, there is a path from $(a,b,c,d,a)$ to $(a,b,d,c,a)$ as well.

  Now we consider $a$-standard colourings as permutations of $\{b,c,d\}$.
  Note that all these permutations can be generated by the transpositions
  $(b,c)$ and $(c,d)$. Therefore, since we can find a path from
  $(a,b,c,d,a)$ to $(a,c,b,d,a)$, and from $(a,b,c,d,a)$ to $(a,b,d,c,a)$,
  there is a path from $(a,b,c,d,a)$ to any other $a$-standard colouring.

  \medskip\noindent
  \textbf{Step 2}:\quad There is a path between any two types of standard
  colourings.\\*[\parskip]
  First, here is a path from an $a$-standard colouring~$\alpha$ to an
  $\alpha(v_3)$-standard colouring\,:

  \begin{center}
    \unitlength0.20mm
    \begin{picture}(120,40)(5,15)
      \put(25,20){\circle*{5}}
      \put(45,20){\circle*{5}}
      \put(65,20){\circle*{5}}
      \put(85,20){\circle*{5}}
      \put(105,20){\circle*{5}}
      \put(25,20){\line(1,0){80}}
      \put(25,33){\makebox(0,0)[b]{$a$}}
      \put(45,33){\makebox(0,0)[b]{$b$}}
      \put(65,33){\makebox(0,0)[b]{$c$}}
      \put(85,33){\makebox(0,0)[b]{$d$}}
      \put(105,33){\makebox(0,0)[b]{$a$}}
    \end{picture}
    \begin{picture}(25,40)(0,15)
      \put(0,20){\vector(1,0){25}}
    \end{picture}
    \begin{picture}(120,40)(5,15)
      \put(25,20){\circle*{5}}
      \put(45,20){\circle*{5}}
      \put(65,20){\circle*{5}}
      \put(85,20){\circle*{5}}
      \put(105,20){\circle*{5}}
      \put(25,20){\line(1,0){80}}
      \put(25,33){\makebox(0,0)[b]{$c$}}
      \put(45,33){\makebox(0,0)[b]{$b$}}
      \put(65,33){\makebox(0,0)[b]{$c$}}
      \put(85,33){\makebox(0,0)[b]{$d$}}
      \put(105,33){\makebox(0,0)[b]{$a$}}
    \end{picture}
    \begin{picture}(25,40)(0,15)
      \put(0,20){\vector(1,0){25}}
    \end{picture}
    \begin{picture}(120,40)(5,15)
      \put(25,20){\circle*{5}}
      \put(45,20){\circle*{5}}
      \put(65,20){\circle*{5}}
      \put(85,20){\circle*{5}}
      \put(105,20){\circle*{5}}
      \put(25,20){\line(1,0){80}}
      \put(25,33){\makebox(0,0)[b]{$c$}}
      \put(45,33){\makebox(0,0)[b]{$b$}}
      \put(65,33){\makebox(0,0)[b]{$a$}}
      \put(85,33){\makebox(0,0)[b]{$d$}}
      \put(105,33){\makebox(0,0)[b]{$a$}}
    \end{picture}
    \begin{picture}(25,40)(0,15)
      \put(0,20){\vector(1,0){25}}
    \end{picture}
    \begin{picture}(120,40)(5,15)
      \put(25,20){\circle*{5}}
      \put(45,20){\circle*{5}}
      \put(65,20){\circle*{5}}
      \put(85,20){\circle*{5}}
      \put(105,20){\circle*{5}}
      \put(25,20){\line(1,0){80}}
      \put(25,33){\makebox(0,0)[b]{$c$}}
      \put(45,33){\makebox(0,0)[b]{$b$}}
      \put(65,33){\makebox(0,0)[b]{$a$}}
      \put(85,33){\makebox(0,0)[b]{$d$}}
      \put(105,33){\makebox(0,0)[b]{$c$}}
    \end{picture}
  \end{center}

  Next, a path from an $a$-standard colouring~$\alpha$ to an
  $\alpha(v_2)$-standard colouring\,:

  \begin{center}
    \unitlength0.20mm
    \begin{picture}(120,40)(5,15)
      \put(25,20){\circle*{5}}
      \put(45,20){\circle*{5}}
      \put(65,20){\circle*{5}}
      \put(85,20){\circle*{5}}
      \put(105,20){\circle*{5}}
      \put(25,20){\line(1,0){80}}
      \put(25,33){\makebox(0,0)[b]{$a$}}
      \put(45,33){\makebox(0,0)[b]{$b$}}
      \put(65,33){\makebox(0,0)[b]{$c$}}
      \put(85,33){\makebox(0,0)[b]{$d$}}
      \put(105,33){\makebox(0,0)[b]{$a$}}
    \end{picture}
    \begin{picture}(25,40)(0,15)
      \put(0,20){\vector(1,0){25}}
    \end{picture}
    \begin{picture}(120,40)(5,15)
      \put(25,20){\circle*{5}}
      \put(45,20){\circle*{5}}
      \put(65,20){\circle*{5}}
      \put(85,20){\circle*{5}}
      \put(105,20){\circle*{5}}
      \put(25,20){\line(1,0){80}}
      \put(25,33){\makebox(0,0)[b]{$d$}}
      \put(45,33){\makebox(0,0)[b]{$b$}}
      \put(65,33){\makebox(0,0)[b]{$c$}}
      \put(85,33){\makebox(0,0)[b]{$d$}}
      \put(105,33){\makebox(0,0)[b]{$a$}}
    \end{picture}
    \begin{picture}(25,40)(0,15)
      \put(0,20){\vector(1,0){25}}
    \end{picture}
    \begin{picture}(120,40)(5,15)
      \put(25,20){\circle*{5}}
      \put(45,20){\circle*{5}}
      \put(65,20){\circle*{5}}
      \put(85,20){\circle*{5}}
      \put(105,20){\circle*{5}}
      \put(25,20){\line(1,0){80}}
      \put(25,33){\makebox(0,0)[b]{$d$}}
      \put(45,33){\makebox(0,0)[b]{$b$}}
      \put(65,33){\makebox(0,0)[b]{$c$}}
      \put(85,33){\makebox(0,0)[b]{$b$}}
      \put(105,33){\makebox(0,0)[b]{$a$}}
    \end{picture}
    \begin{picture}(25,40)(0,15)
      \put(0,20){\vector(1,0){25}}
    \end{picture}
    \begin{picture}(120,40)(5,15)
      \put(25,20){\circle*{5}}
      \put(45,20){\circle*{5}}
      \put(65,20){\circle*{5}}
      \put(85,20){\circle*{5}}
      \put(105,20){\circle*{5}}
      \put(25,20){\line(1,0){80}}
      \put(25,33){\makebox(0,0)[b]{$d$}}
      \put(45,33){\makebox(0,0)[b]{$a$}}
      \put(65,33){\makebox(0,0)[b]{$c$}}
      \put(85,33){\makebox(0,0)[b]{$b$}}
      \put(105,33){\makebox(0,0)[b]{$a$}}
    \end{picture}
    \begin{picture}(25,40)(0,15)
      \put(0,20){\vector(1,0){25}}
    \end{picture}
    \begin{picture}(120,40)(5,15)
      \put(25,20){\circle*{5}}
      \put(45,20){\circle*{5}}
      \put(65,20){\circle*{5}}
      \put(85,20){\circle*{5}}
      \put(105,20){\circle*{5}}
      \put(25,20){\line(1,0){80}}
      \put(25,33){\makebox(0,0)[b]{$d$}}
      \put(45,33){\makebox(0,0)[b]{$a$}}
      \put(65,33){\makebox(0,0)[b]{$c$}}
      \put(85,33){\makebox(0,0)[b]{$b$}}
      \put(105,33){\makebox(0,0)[b]{$d$}}
    \end{picture}
    \begin{picture}(25,40)(0,15)
      \put(0,20){\vector(1,0){25}}
    \end{picture}
    \begin{picture}(120,40)(5,15)
      \put(25,20){\circle*{5}}
      \put(45,20){\circle*{5}}
      \put(65,20){\circle*{5}}
      \put(85,20){\circle*{5}}
      \put(105,20){\circle*{5}}
      \put(25,20){\line(1,0){80}}
      \put(25,33){\makebox(0,0)[b]{$b$}}
      \put(45,33){\makebox(0,0)[b]{$a$}}
      \put(65,33){\makebox(0,0)[b]{$c$}}
      \put(85,33){\makebox(0,0)[b]{$b$}}
      \put(105,33){\makebox(0,0)[b]{$d$}}
    \end{picture}
    \begin{picture}(25,40)(0,15)
      \put(0,20){\vector(1,0){25}}
    \end{picture}
    \begin{picture}(120,40)(5,15)
      \put(25,20){\circle*{5}}
      \put(45,20){\circle*{5}}
      \put(65,20){\circle*{5}}
      \put(85,20){\circle*{5}}
      \put(105,20){\circle*{5}}
      \put(25,20){\line(1,0){80}}
      \put(25,33){\makebox(0,0)[b]{$b$}}
      \put(45,33){\makebox(0,0)[b]{$a$}}
      \put(65,33){\makebox(0,0)[b]{$c$}}
      \put(85,33){\makebox(0,0)[b]{$a$}}
      \put(105,33){\makebox(0,0)[b]{$d$}}
    \end{picture}
    \begin{picture}(25,40)(0,15)
      \put(0,20){\vector(1,0){25}}
    \end{picture}
    \begin{picture}(120,40)(5,15)
      \put(25,20){\circle*{5}}
      \put(45,20){\circle*{5}}
      \put(65,20){\circle*{5}}
      \put(85,20){\circle*{5}}
      \put(105,20){\circle*{5}}
      \put(25,20){\line(1,0){80}}
      \put(25,33){\makebox(0,0)[b]{$b$}}
      \put(45,33){\makebox(0,0)[b]{$d$}}
      \put(65,33){\makebox(0,0)[b]{$c$}}
      \put(85,33){\makebox(0,0)[b]{$a$}}
      \put(105,33){\makebox(0,0)[b]{$d$}}
    \end{picture}
    \begin{picture}(25,40)(0,15)
      \put(0,20){\vector(1,0){25}}
    \end{picture}
    \begin{picture}(120,40)(5,15)
      \put(25,20){\circle*{5}}
      \put(45,20){\circle*{5}}
      \put(65,20){\circle*{5}}
      \put(85,20){\circle*{5}}
      \put(105,20){\circle*{5}}
      \put(25,20){\line(1,0){80}}
      \put(25,33){\makebox(0,0)[b]{$b$}}
      \put(45,33){\makebox(0,0)[b]{$d$}}
      \put(65,33){\makebox(0,0)[b]{$c$}}
      \put(85,33){\makebox(0,0)[b]{$a$}}
      \put(105,33){\makebox(0,0)[b]{$b$}}
    \end{picture}
  \end{center}

  By symmetry, there is also a path from an $a$-standard colouring~$\alpha$
  to an $\alpha(v_4)$-standard colouring.

  \medskip\noindent
  \textbf{Step 3}:\quad Each colouring has a path to some standard
  colouring.\\*[\parskip]
  Let~$\alpha$ be a strong 4-colouring of~$P_5$. Then~$\alpha$ has one of
  the following forms\,:

  \begin{center}
    \unitlength0.20mm
    \begin{picture}(120,40)(5,15)
      \put(25,20){\circle*{5}}
      \put(45,20){\circle*{5}}
      \put(65,20){\circle*{5}}
      \put(85,20){\circle*{5}}
      \put(105,20){\circle*{5}}
      \put(25,20){\line(1,0){80}}
      \put(25,33){\makebox(0,0)[b]{$a$}}
      \put(45,33){\makebox(0,0)[b]{$b$}}
      \put(65,33){\makebox(0,0)[b]{$c$}}
      \put(85,33){\makebox(0,0)[b]{$d$}}
      \put(105,33){\makebox(0,0)[b]{$a$}}
    \end{picture}
    \begin{picture}(120,40)(5,15)
      \put(25,20){\circle*{5}}
      \put(45,20){\circle*{5}}
      \put(65,20){\circle*{5}}
      \put(85,20){\circle*{5}}
      \put(105,20){\circle*{5}}
      \put(25,20){\line(1,0){80}}
      \put(25,33){\makebox(0,0)[b]{$a$}}
      \put(45,33){\makebox(0,0)[b]{$b$}}
      \put(65,33){\makebox(0,0)[b]{$c$}}
      \put(85,33){\makebox(0,0)[b]{$a$}}
      \put(105,33){\makebox(0,0)[b]{$d$}}
    \end{picture}
    \begin{picture}(120,40)(5,15)
      \put(25,20){\circle*{5}}
      \put(45,20){\circle*{5}}
      \put(65,20){\circle*{5}}
      \put(85,20){\circle*{5}}
      \put(105,20){\circle*{5}}
      \put(25,20){\line(1,0){80}}
      \put(25,33){\makebox(0,0)[b]{$a$}}
      \put(45,33){\makebox(0,0)[b]{$b$}}
      \put(65,33){\makebox(0,0)[b]{$a$}}
      \put(85,33){\makebox(0,0)[b]{$c$}}
      \put(105,33){\makebox(0,0)[b]{$d$}}
    \end{picture}
    \begin{picture}(120,40)(5,15)
      \put(25,20){\circle*{5}}
      \put(45,20){\circle*{5}}
      \put(65,20){\circle*{5}}
      \put(85,20){\circle*{5}}
      \put(105,20){\circle*{5}}
      \put(25,20){\line(1,0){80}}
      \put(25,33){\makebox(0,0)[b]{$b$}}
      \put(45,33){\makebox(0,0)[b]{$a$}}
      \put(65,33){\makebox(0,0)[b]{$c$}}
      \put(85,33){\makebox(0,0)[b]{$d$}}
      \put(105,33){\makebox(0,0)[b]{$a$}}
    \end{picture}
    \begin{picture}(120,40)(5,15)
      \put(25,20){\circle*{5}}
      \put(45,20){\circle*{5}}
      \put(65,20){\circle*{5}}
      \put(85,20){\circle*{5}}
      \put(105,20){\circle*{5}}
      \put(25,20){\line(1,0){80}}
      \put(25,33){\makebox(0,0)[b]{$b$}}
      \put(45,33){\makebox(0,0)[b]{$a$}}
      \put(65,33){\makebox(0,0)[b]{$c$}}
      \put(85,33){\makebox(0,0)[b]{$a$}}
      \put(105,33){\makebox(0,0)[b]{$d$}}
    \end{picture}
    \begin{picture}(120,40)(5,15)
      \put(25,20){\circle*{5}}
      \put(45,20){\circle*{5}}
      \put(65,20){\circle*{5}}
      \put(85,20){\circle*{5}}
      \put(105,20){\circle*{5}}
      \put(25,20){\line(1,0){80}}
      \put(25,33){\makebox(0,0)[b]{$b$}}
      \put(45,33){\makebox(0,0)[b]{$c$}}
      \put(65,33){\makebox(0,0)[b]{$a$}}
      \put(85,33){\makebox(0,0)[b]{$d$}}
      \put(105,33){\makebox(0,0)[b]{$a$}}
    \end{picture}
  \end{center}

  The first form already is an $a$-standard colouring; for the second and
  the third ones we just recolour the vertex~$v_1$ to~$d$; while for the
  fourth and the sixth ones we just recolour the vertex~$v_5$ to~$b$.
  Finally, for the fifth form the following is a path to a $b$-standard
  colouring\,:

  \begin{center}
   \unitlength0.20mm
    \begin{picture}(120,40)(5,15)
      \put(25,20){\circle*{5}}
      \put(45,20){\circle*{5}}
      \put(65,20){\circle*{5}}
      \put(85,20){\circle*{5}}
      \put(105,20){\circle*{5}}
      \put(25,20){\line(1,0){80}}
      \put(25,33){\makebox(0,0)[b]{$b$}}
      \put(45,33){\makebox(0,0)[b]{$a$}}
      \put(65,33){\makebox(0,0)[b]{$c$}}
      \put(85,33){\makebox(0,0)[b]{$a$}}
      \put(105,33){\makebox(0,0)[b]{$d$}}
    \end{picture}
    \begin{picture}(25,40)(0,15)
      \put(0,20){\vector(1,0){25}}
    \end{picture}
    \begin{picture}(120,40)(5,15)
      \put(25,20){\circle*{5}}
      \put(45,20){\circle*{5}}
      \put(65,20){\circle*{5}}
      \put(85,20){\circle*{5}}
      \put(105,20){\circle*{5}}
      \put(25,20){\line(1,0){80}}
      \put(25,33){\makebox(0,0)[b]{$b$}}
      \put(45,33){\makebox(0,0)[b]{$d$}}
      \put(65,33){\makebox(0,0)[b]{$c$}}
      \put(85,33){\makebox(0,0)[b]{$a$}}
      \put(105,33){\makebox(0,0)[b]{$d$}}
    \end{picture}
    \begin{picture}(25,40)(0,15)
      \put(0,20){\vector(1,0){25}}
    \end{picture}
    \begin{picture}(120,40)(5,15)
      \put(25,20){\circle*{5}}
      \put(45,20){\circle*{5}}
      \put(65,20){\circle*{5}}
      \put(85,20){\circle*{5}}
      \put(105,20){\circle*{5}}
      \put(25,20){\line(1,0){80}}
      \put(25,33){\makebox(0,0)[b]{$b$}}
      \put(45,33){\makebox(0,0)[b]{$d$}}
      \put(65,33){\makebox(0,0)[b]{$c$}}
      \put(85,33){\makebox(0,0)[b]{$a$}}
      \put(105,33){\makebox(0,0)[b]{$b$}}
    \end{picture}
  \end{center}

  \medskip\noindent
  It is straightforward to see that appropriate renaming of the colours and
  sequence of the paths in Steps 1\,--\,3 will transform any strong
  4-colouring of~$P_5$ into any other strong 4-colouring.
\end{proof}

\noindent
  We extend the last result by showing that $S_k(P_{k+1})$ is connected, for
  all $k\ge4$.

\begin{prop} \label{B.9}
  For all $k\ge4$, $S_k(P_{k+1})$ is connected.
\end{prop}

\begin{proof}
  We will prove this by induction on~$k$. We have already shown the
  proposition is true for $k=4$.

  Let~$\alpha$ and~$\beta$ be strong $k$-colourings of
  $P_{k+1}=v_1v_2...v_{k+1}$, for some $k\ge5$. We can assume that
  in~$\alpha$, the vertex $v_{k+1}$ has a unique colour. Otherwise, there
  is another vertex~$v_i$ such that $\alpha(v_i)=\alpha(v_{k+1})$. Then
  just recolour~$v_i$ to a colour different from $\alpha(v_{k+1})$. Next,
  we can assume that this unique colour on $v_{k+1}$ in~$\alpha$ is~$a$.

  If the vertex $v_{k+1}$ is the only vertex coloured~$a$ in~$\beta$ as
  well, then we can just remove $v_{k+1}$. Let~$\alpha'$ and~$\beta'$ be
  the strong $k$-colourings of $P_k=v_1v_2...v_k$ obtained from~$\alpha$
  and~$\beta$, respectively. Since, by induction, $S_{k-1}(P_k)$ is
  connected, there is a path from~$\alpha'$ to~$\beta'$ in $S_{k-1}(P_k)$.
  Using the same steps on $P_{k+1}$ gives a path from~$\alpha$ to~$\beta$
  in $S_k(P_{k+1})$.

  So we can assume that in~$\beta$, $v_{k+1}$ is not coloured~$a$ or is not
  the only vertex coloured~$a$. We distinguish~4 cases.

  \medskip\noindent
  \textbf{Case 1}:\quad In~$\beta$, $v_{k+1}$ is coloured~$a$, but there is
  a second vertex~$v_i$ coloured~$a$ as well.\\*[\parskip]
  Then just recolour~$v_i$ to some different from~$a$, and so $v_{k+1}$ is
  now the only vertex coloured~$a$. We are done by the paragraph above.

  \medskip\noindent
  \textbf{Case 2}:\quad In~$\beta$, $v_{k+1}$ and some other vertex~$v_i$
  have the same colour $b\ne a$, while a third vertex~$v_j$ is
  coloured~$a$.

  \medskip\noindent
  \textbf{Subcase 2.1}:\quad $v_{k+1}$ and~$v_j$ are not
  adjacent.\\*[\parskip]
  Then just recolour $v_{k+1}$ to~$a$, and we are back to Case~1.

  \medskip\noindent
  \textbf{Subcase 2.2}:\quad $v_{k+1}$ and~$v_j$ are adjacent, i.e.,
  $j=k$.\\*[\parskip]
  Call a colouring of $S_k(P_{k+1})$ \emph{good} if we can recolour~$v_i$
  to a colour which is not one of
  $\{\beta(v_{k-1}),\beta(v_k)=a,\beta(v_{k+1})=b\}$.

  Note that~$\beta$ is good when $k\ge6$, or $k=5$ and $i\ne2$. If~$\beta$
  is good, we can recolour~$v_i$ to the colour $\beta(v_l)$ for some
  $l\notin\{i-1,i,i+1,k-1,k,k+1\}$, to obtain the strong
  $k$-colouring~$\gamma$, Let~$\delta$ be the strong $k$-colouring of
  $P_{k+1}$, obtained from~$\gamma$ by swapping the colours of~$v_i$
  and~$v_k$. By ignoring the vertex $v_{k+1}$, we can consider~$\gamma$
  and~$\delta$ as strong $(k-1)$-colourings of~$P_k$. Since $S_{k-1}(P_k)$
  is connected, there is a path between these two colourings. We then apply
  this path to a path in $S_k(P_{k+1})$ from~$\gamma$ to~$\delta$.

  Next, we will form a path from~$\delta$ to a colouring in which vertex
  $v_{k+1}$ is the only vertex coloured~$a$. Therefore, we also have a path
  from~$\beta$ to this colouring. In~$\delta$, we first recolour~$v_l$ to
  $\beta(v_{k+1})=b$ and then recolour $v_{k+1}$ to~$a$. Finally, recolour
  vertex~$v_i$, which is previously coloured~$a$, to another colour, and we
  are done.

  We now suppose that~$\beta$ is not good, i.e., $k=5$ and $i=2$. Then
  there is a path from~$\beta$ to a colouring in which $v_{k+1}$ is the
  only vertex coloured~$a$.

  \begin{center}
    \unitlength0.20mm
    \begin{picture}(140,40)(5,15)
      \put(25,20){\circle*{5}}
      \put(45,20){\circle*{5}}
      \put(65,20){\circle*{5}}
      \put(85,20){\circle*{5}}
      \put(105,20){\circle*{5}}
      \put(125,20){\circle*{5}}
      \put(25,20){\line(1,0){100}}
      \put(25,35){\makebox(0,0)[b]{$c$}}
      \put(45,35){\makebox(0,0)[b]{$b$}}
      \put(65,35){\makebox(0,0)[b]{$d$}}
      \put(85,35){\makebox(0,0)[b]{$e$}}
      \put(105,35){\makebox(0,0)[b]{$a$}}
      \put(125,35){\makebox(0,0)[b]{$b$}}
    \end{picture}
    \begin{picture}(25,40)(0,15)
      \put(0,20){\vector(1,0){25}}
    \end{picture}
    \begin{picture}(140,40)(5,15)
      \put(25,20){\circle*{5}}
      \put(45,20){\circle*{5}}
      \put(65,20){\circle*{5}}
      \put(85,20){\circle*{5}}
      \put(105,20){\circle*{5}}
      \put(125,20){\circle*{5}}
      \put(25,20){\line(1,0){100}}
      \put(25,35){\makebox(0,0)[b]{$c$}}
      \put(45,35){\makebox(0,0)[b]{$a$}}
      \put(65,35){\makebox(0,0)[b]{$d$}}
      \put(85,35){\makebox(0,0)[b]{$e$}}
      \put(105,35){\makebox(0,0)[b]{$a$}}
      \put(125,35){\makebox(0,0)[b]{$b$}}
    \end{picture}
    \begin{picture}(25,40)(0,15)
      \put(0,20){\vector(1,0){25}}
    \end{picture}
    \begin{picture}(140,40)(5,15)
      \put(25,20){\circle*{5}}
      \put(45,20){\circle*{5}}
      \put(65,20){\circle*{5}}
      \put(85,20){\circle*{5}}
      \put(105,20){\circle*{5}}
      \put(125,20){\circle*{5}}
      \put(25,20){\line(1,0){100}}
      \put(25,35){\makebox(0,0)[b]{$c$}}
      \put(45,35){\makebox(0,0)[b]{$a$}}
      \put(65,35){\makebox(0,0)[b]{$d$}}
      \put(85,35){\makebox(0,0)[b]{$e$}}
      \put(105,35){\makebox(0,0)[b]{$d$}}
      \put(125,35){\makebox(0,0)[b]{$b$}}
    \end{picture}
    \begin{picture}(25,40)(0,15)
      \put(0,20){\vector(1,0){25}}
    \end{picture}
    \begin{picture}(140,40)(5,15)
      \put(25,20){\circle*{5}}
      \put(45,20){\circle*{5}}
      \put(65,20){\circle*{5}}
      \put(85,20){\circle*{5}}
      \put(105,20){\circle*{5}}
      \put(125,20){\circle*{5}}
      \put(25,20){\line(1,0){100}}
      \put(25,35){\makebox(0,0)[b]{$c$}}
      \put(45,35){\makebox(0,0)[b]{$a$}}
      \put(65,35){\makebox(0,0)[b]{$b$}}
      \put(85,35){\makebox(0,0)[b]{$e$}}
      \put(105,35){\makebox(0,0)[b]{$d$}}
      \put(125,35){\makebox(0,0)[b]{$b$}}
    \end{picture}
    \begin{picture}(25,40)(0,15)
      \put(0,20){\vector(1,0){25}}
    \end{picture}
    \begin{picture}(140,40)(5,15)
      \put(25,20){\circle*{5}}
      \put(45,20){\circle*{5}}
      \put(65,20){\circle*{5}}
      \put(85,20){\circle*{5}}
      \put(105,20){\circle*{5}}
      \put(125,20){\circle*{5}}
      \put(25,20){\line(1,0){100}}
      \put(25,35){\makebox(0,0)[b]{$c$}}
      \put(45,35){\makebox(0,0)[b]{$a$}}
      \put(65,35){\makebox(0,0)[b]{$b$}}
      \put(85,35){\makebox(0,0)[b]{$e$}}
      \put(105,35){\makebox(0,0)[b]{$d$}}
      \put(125,35){\makebox(0,0)[b]{$a$}}
    \end{picture}
    \begin{picture}(25,40)(0,15)
      \put(0,20){\vector(1,0){25}}
    \end{picture}
    \begin{picture}(140,40)(5,15)
      \put(25,20){\circle*{5}}
      \put(45,20){\circle*{5}}
      \put(65,20){\circle*{5}}
      \put(85,20){\circle*{5}}
      \put(105,20){\circle*{5}}
      \put(125,20){\circle*{5}}
      \put(25,20){\line(1,0){100}}
      \put(25,35){\makebox(0,0)[b]{$c$}}
      \put(45,35){\makebox(0,0)[b]{$e$}}
      \put(65,35){\makebox(0,0)[b]{$b$}}
      \put(85,35){\makebox(0,0)[b]{$e$}}
      \put(105,35){\makebox(0,0)[b]{$d$}}
      \put(125,35){\makebox(0,0)[b]{$a$}}
    \end{picture}
  \end{center}

  \medskip\noindent
  \textbf{Case 3}:\quad In~$\beta$, $v_i$ and~$v_j$ are coloured~$a$ for
  some $i,j\ne k+1$.\\*[\parskip]
  Without loss of generality, we may assume that~$v_i$ is not adjacent to
  $v_{k+1}$. Then we recolour~$v_i$ to $\beta(v_{k+1})$, and we are back to
  Case~2.

  \medskip\noindent
  \textbf{Case 4}:\quad In~$\beta$, $v_i$ and~$v_j$ have the same colour
  $b\ne a$, a third vertex~$v_\ell$ is coloured~$a$, for some
  $i,j,\ell\ne k+1$.\\*[\parskip]
  Without loss of generality, we may assume that~$v_i$ is not adjacent to
  $v_{k+1}$. Then we recolour~$v_i$ to $\beta(v_{k+1})$, and we are back to
  Case~2.
\end{proof}

\noindent
  Combining it all, we get the promised result on the strong colour graph of
  paths.

\begin{thm}\label{B.10}
  The strong colour graph $S_k(P_n)$ is connected if and only if $k\ge3$,
  $n\ge5$ and $n\ge k+1$.
\end{thm}

\begin{proof}
  We already have seen that $S_3(P_4)$ and $S_2(P_n)$, $n\ge3$, are not
  connected, while $S_3(P_5)$ and $S_k(P_{k+1})$, $k\ge4$, are connected.
  Applying Theorem~\ref{Z.6} completes the proof.
\end{proof}

\section{The Strong $k$-Colour Graph of Cycles}\label{sec4}

  In this section we want to show that the strong $k$-colour graph of a cycle
  with~$n$ vertices, $S_k(C_n)$, is connected if and only if $k\ge4$, $n\ge6$
  and $n\ge k+1$. Before we prove the theorem, we prove some tools used in
  this proof.

  To orient a cycle means to orient each edge on the cycle so that a directed
  cycle is obtained. If~$C$ is a cycle, then by~$\arrc$ we denote the cycle
  with one of the two possible orientations of~$d$. Given a
  3-colouring~$\alpha$ using colours $\{1,2,3\}$, the weight of an
  edge~$e=uv$ oriented from~$u$ to~$v$ is
  \[w(\overrightarrow{uv},\alpha)\:=\:\left\{\begin{array}{ll}
    +1,&\text{if $\alpha(u)\alpha(v)\in\{12,\,23,\,31\}$};\\
    -1,&\text{if $\alpha(u)\alpha(v)\in\{21,\,32,\,13\}$}.\\
  \end{array}\right.\]
  The weight $W(\arrc,\alpha)$ of an oriented cycle~$\arrc$ is the sum of the
  weights of its oriented edges.

\begin{lem}\label{Lemma5Jan}
  \textsc{(Cereceda et al.~\cite{vertex-colouring})} Let~$\alpha$ be a
  3-colouring of a graph~$G$ that contains a cycle~$C$. If
  $W(\arrc,\alpha)\ne0$, then $C_k(G)$ is not connected.
\end{lem}

\begin{prop}\label{A.1.5}
  For all $n\ge3$, $S_3(C_n)$ is not connected.
\end{prop}

\begin{proof}
  By Lemmas~\ref{Z.10} and~\ref{Lemma5Jan}, it is enough to find a strong
  3-colouring~$\alpha$ with $W(\overrightarrow{C_n},\alpha)\ne0$. If
  $n=3\,\ell$ for some positive integer~$\ell$, the pattern
  1,2,3,1,2,3,...,1,2,3 provides a 3-colouring~$\alpha$ of~$C_n$ with
  $W(\overrightarrow{C_n},\alpha)=n\ne0$. For $n=4$, it is easy to see that
  $S_3(C_4)$ is a graph with 12 isolated vertices. If $n=3\,\ell+1>4$, then
  we use the pattern 1,2,3,1,2,3,...,1,2,3,2, which gives
  $W(\overrightarrow{C_n},\alpha)=n-4\ne0$. Finally, if $n=3\,\ell+2\ge5$,
  then we use the pattern 1,2,3,1,2,3,...,1,2,3,1,2, with
  $W(\overrightarrow{C_n},\alpha)=n-2\ne0$.
\end{proof}

\begin{prop}\label{A.2.2}
  The strong colour graph $S_4(C_5)$ is not connected.
\end{prop}

\begin{proof}
  For any strong 4-colouring of the 5-cycle~$C_5$, there are only two
  vertices having the same colour. Thus each strong 4-vertex-colouring
  of~$C_5$ can be recoloured only on these two vertices, and each of these
  two vertices can be recoloured to only one new colour (\,since the two
  different colours of their neighbours are forbidden\,). This means each
  colouring has degree two in $S_4(C_5)$.

  Straightforward counting shows that $S_4(C_5)$ has 120 vertices. But each
  colouring in $S_4(C_5)$ is contained in some cycle of length~20. To see
  this, we start with some strong 4-colouring of~$C_5$ and recolour\,:

  \begin{center}
    \unitlength0.20mm
    \begin{picture}(120,50)(5,5)
      \put(25,20){\circle*{5}}
      \put(45,20){\circle*{5}}
      \put(65,20){\circle*{5}}
      \put(85,20){\circle*{5}}
      \put(105,20){\circle*{5}}
      \put(25,20){\line(1,0){80}}
      \qbezier(25,20)(65,-5)(105,20)
      \put(25,33){\makebox(0,0)[b]{$a$}}
      \put(45,33){\makebox(0,0)[b]{$b$}}
      \put(65,33){\makebox(0,0)[b]{$a$}}
      \put(85,33){\makebox(0,0)[b]{$c$}}
      \put(105,33){\makebox(0,0)[b]{$d$}}
    \end{picture}
    \begin{picture}(25,50)(0,5)
      \put(0,20){\vector(1,0){25}}
    \end{picture}
    \begin{picture}(120,50)(5,5)
      \put(25,20){\circle*{5}}
      \put(45,20){\circle*{5}}
      \put(65,20){\circle*{5}}
      \put(85,20){\circle*{5}}
      \put(105,20){\circle*{5}}
      \put(25,20){\line(1,0){80}}
      \qbezier(25,20)(65,-5)(105,20)
      \put(25,33){\makebox(0,0)[b]{$a$}}
      \put(45,33){\makebox(0,0)[b]{$b$}}
      \put(65,33){\makebox(0,0)[b]{$d$}}
      \put(85,33){\makebox(0,0)[b]{$c$}}
      \put(105,33){\makebox(0,0)[b]{$d$}}
    \end{picture}
    \begin{picture}(25,50)(0,5)
      \put(0,20){\vector(1,0){25}}
    \end{picture}
    \begin{picture}(120,50)(5,5)
      \put(25,20){\circle*{5}}
      \put(45,20){\circle*{5}}
      \put(65,20){\circle*{5}}
      \put(85,20){\circle*{5}}
      \put(105,20){\circle*{5}}
      \put(25,20){\line(1,0){80}}
      \qbezier(25,20)(65,-5)(105,20)
      \put(25,33){\makebox(0,0)[b]{$a$}}
      \put(45,33){\makebox(0,0)[b]{$b$}}
      \put(65,33){\makebox(0,0)[b]{$d$}}
      \put(85,33){\makebox(0,0)[b]{$c$}}
      \put(105,33){\makebox(0,0)[b]{$b$}}
    \end{picture}
    \begin{picture}(25,50)(0,5)
      \put(0,20){\vector(1,0){25}}
    \end{picture}
    \begin{picture}(120,50)(5,5)
      \put(25,20){\circle*{5}}
      \put(45,20){\circle*{5}}
      \put(65,20){\circle*{5}}
      \put(85,20){\circle*{5}}
      \put(105,20){\circle*{5}}
      \put(25,20){\line(1,0){80}}
      \qbezier(25,20)(65,-5)(105,20)
      \put(25,33){\makebox(0,0)[b]{$a$}}
      \put(45,33){\makebox(0,0)[b]{$c$}}
      \put(65,33){\makebox(0,0)[b]{$d$}}
      \put(85,33){\makebox(0,0)[b]{$c$}}
      \put(105,33){\makebox(0,0)[b]{$b$}}
    \end{picture}
    \begin{picture}(25,50)(0,5)
      \put(0,20){\vector(1,0){25}}
    \end{picture}
    \begin{picture}(120,50)(5,5)
      \put(25,20){\circle*{5}}
      \put(45,20){\circle*{5}}
      \put(65,20){\circle*{5}}
      \put(85,20){\circle*{5}}
      \put(105,20){\circle*{5}}
      \put(25,20){\line(1,0){80}}
      \qbezier(25,20)(65,-5)(105,20)
      \put(25,33){\makebox(0,0)[b]{$a$}}
      \put(45,33){\makebox(0,0)[b]{$c$}}
      \put(65,33){\makebox(0,0)[b]{$d$}}
      \put(85,33){\makebox(0,0)[b]{$a$}}
      \put(105,33){\makebox(0,0)[b]{$b$}}
    \end{picture}
    \begin{picture}(25,50)(0,5)
      \put(0,20){\vector(1,0){25}}
    \end{picture}
    \begin{picture}(120,50)(5,5)
      \put(25,20){\circle*{5}}
      \put(45,20){\circle*{5}}
      \put(65,20){\circle*{5}}
      \put(85,20){\circle*{5}}
      \put(105,20){\circle*{5}}
      \put(25,20){\line(1,0){80}}
      \qbezier(25,20)(65,-5)(105,20)
      \put(25,33){\makebox(0,0)[b]{$d$}}
      \put(45,33){\makebox(0,0)[b]{$c$}}
      \put(65,33){\makebox(0,0)[b]{$d$}}
      \put(85,33){\makebox(0,0)[b]{$a$}}
      \put(105,33){\makebox(0,0)[b]{$b$}}
    \end{picture}
    \begin{picture}(25,50)(0,5)
      \put(0,20){\vector(1,0){25}}
    \end{picture}
    \begin{picture}(120,50)(5,5)
      \put(25,20){\circle*{5}}
      \put(45,20){\circle*{5}}
      \put(65,20){\circle*{5}}
      \put(85,20){\circle*{5}}
      \put(105,20){\circle*{5}}
      \put(25,20){\line(1,0){80}}
      \qbezier(25,20)(65,-5)(105,20)
      \put(25,33){\makebox(0,0)[b]{$d$}}
      \put(45,33){\makebox(0,0)[b]{$c$}}
      \put(65,33){\makebox(0,0)[b]{$b$}}
      \put(85,33){\makebox(0,0)[b]{$a$}}
      \put(105,33){\makebox(0,0)[b]{$b$}}
    \end{picture}
    \begin{picture}(25,50)(0,5)
      \put(0,20){\vector(1,0){25}}
    \end{picture}
    \begin{picture}(120,50)(5,5)
      \put(25,20){\circle*{5}}
      \put(45,20){\circle*{5}}
      \put(65,20){\circle*{5}}
      \put(85,20){\circle*{5}}
      \put(105,20){\circle*{5}}
      \put(25,20){\line(1,0){80}}
      \qbezier(25,20)(65,-5)(105,20)
      \put(25,33){\makebox(0,0)[b]{$d$}}
      \put(45,33){\makebox(0,0)[b]{$c$}}
      \put(65,33){\makebox(0,0)[b]{$b$}}
      \put(85,33){\makebox(0,0)[b]{$a$}}
      \put(105,33){\makebox(0,0)[b]{$c$}}
    \end{picture}
    \begin{picture}(25,50)(0,5)
      \put(0,20){\vector(1,0){25}}
    \end{picture}
    \begin{picture}(120,50)(5,5)
      \put(25,20){\circle*{5}}
      \put(45,20){\circle*{5}}
      \put(65,20){\circle*{5}}
      \put(85,20){\circle*{5}}
      \put(105,20){\circle*{5}}
      \put(25,20){\line(1,0){80}}
      \qbezier(25,20)(65,-5)(105,20)
      \put(25,33){\makebox(0,0)[b]{$d$}}
      \put(45,33){\makebox(0,0)[b]{$a$}}
      \put(65,33){\makebox(0,0)[b]{$b$}}
      \put(85,33){\makebox(0,0)[b]{$a$}}
      \put(105,33){\makebox(0,0)[b]{$c$}}
    \end{picture}
    \begin{picture}(25,50)(0,5)
      \put(0,20){\vector(1,0){25}}
    \end{picture}
    \begin{picture}(120,50)(5,5)
      \put(25,20){\circle*{5}}
      \put(45,20){\circle*{5}}
      \put(65,20){\circle*{5}}
      \put(85,20){\circle*{5}}
      \put(105,20){\circle*{5}}
      \put(25,20){\line(1,0){80}}
      \qbezier(25,20)(65,-5)(105,20)
      \put(25,33){\makebox(0,0)[b]{$d$}}
      \put(45,33){\makebox(0,0)[b]{$a$}}
      \put(65,33){\makebox(0,0)[b]{$b$}}
      \put(85,33){\makebox(0,0)[b]{$d$}}
      \put(105,33){\makebox(0,0)[b]{$c$}}
    \end{picture}
    \begin{picture}(25,50)(0,5)
      \put(0,20){\vector(1,0){25}}
    \end{picture}
    \begin{picture}(120,50)(5,5)
      \put(25,20){\circle*{5}}
      \put(45,20){\circle*{5}}
      \put(65,20){\circle*{5}}
      \put(85,20){\circle*{5}}
      \put(105,20){\circle*{5}}
      \put(25,20){\line(1,0){80}}
      \qbezier(25,20)(65,-5)(105,20)
      \put(25,33){\makebox(0,0)[b]{$b$}}
      \put(45,33){\makebox(0,0)[b]{$a$}}
      \put(65,33){\makebox(0,0)[b]{$b$}}
      \put(85,33){\makebox(0,0)[b]{$d$}}
      \put(105,33){\makebox(0,0)[b]{$c$}}
    \end{picture}
    \begin{picture}(25,50)(0,5)
      \put(0,20){\vector(1,0){25}}
    \end{picture}
    \begin{picture}(120,50)(5,5)
      \put(25,20){\circle*{5}}
      \put(45,20){\circle*{5}}
      \put(65,20){\circle*{5}}
      \put(85,20){\circle*{5}}
      \put(105,20){\circle*{5}}
      \put(25,20){\line(1,0){80}}
      \qbezier(25,20)(65,-5)(105,20)
      \put(25,33){\makebox(0,0)[b]{$b$}}
      \put(45,33){\makebox(0,0)[b]{$a$}}
      \put(65,33){\makebox(0,0)[b]{$c$}}
      \put(85,33){\makebox(0,0)[b]{$d$}}
      \put(105,33){\makebox(0,0)[b]{$c$}}
    \end{picture}
    \begin{picture}(25,50)(0,5)
      \put(0,20){\vector(1,0){25}}
    \end{picture}
    \begin{picture}(120,50)(5,5)
      \put(25,20){\circle*{5}}
      \put(45,20){\circle*{5}}
      \put(65,20){\circle*{5}}
      \put(85,20){\circle*{5}}
      \put(105,20){\circle*{5}}
      \put(25,20){\line(1,0){80}}
      \qbezier(25,20)(65,-5)(105,20)
      \put(25,33){\makebox(0,0)[b]{$b$}}
      \put(45,33){\makebox(0,0)[b]{$a$}}
      \put(65,33){\makebox(0,0)[b]{$c$}}
      \put(85,33){\makebox(0,0)[b]{$d$}}
      \put(105,33){\makebox(0,0)[b]{$a$}}
    \end{picture}
    \begin{picture}(25,50)(0,5)
      \put(0,20){\vector(1,0){25}}
    \end{picture}
    \begin{picture}(120,50)(5,5)
      \put(25,20){\circle*{5}}
      \put(45,20){\circle*{5}}
      \put(65,20){\circle*{5}}
      \put(85,20){\circle*{5}}
      \put(105,20){\circle*{5}}
      \put(25,20){\line(1,0){80}}
      \qbezier(25,20)(65,-5)(105,20)
      \put(25,33){\makebox(0,0)[b]{$b$}}
      \put(45,33){\makebox(0,0)[b]{$d$}}
      \put(65,33){\makebox(0,0)[b]{$c$}}
      \put(85,33){\makebox(0,0)[b]{$d$}}
      \put(105,33){\makebox(0,0)[b]{$a$}}
    \end{picture}
    \begin{picture}(25,50)(0,5)
      \put(0,20){\vector(1,0){25}}
    \end{picture}
    \begin{picture}(120,50)(5,5)
      \put(25,20){\circle*{5}}
      \put(45,20){\circle*{5}}
      \put(65,20){\circle*{5}}
      \put(85,20){\circle*{5}}
      \put(105,20){\circle*{5}}
      \put(25,20){\line(1,0){80}}
      \qbezier(25,20)(65,-5)(105,20)
      \put(25,33){\makebox(0,0)[b]{$b$}}
      \put(45,33){\makebox(0,0)[b]{$d$}}
      \put(65,33){\makebox(0,0)[b]{$c$}}
      \put(85,33){\makebox(0,0)[b]{$b$}}
      \put(105,33){\makebox(0,0)[b]{$a$}}
    \end{picture}
    \begin{picture}(25,50)(0,5)
      \put(0,20){\vector(1,0){25}}
    \end{picture}
    \begin{picture}(120,50)(5,5)
      \put(25,20){\circle*{5}}
      \put(45,20){\circle*{5}}
      \put(65,20){\circle*{5}}
      \put(85,20){\circle*{5}}
      \put(105,20){\circle*{5}}
      \put(25,20){\line(1,0){80}}
      \qbezier(25,20)(65,-5)(105,20)
      \put(25,33){\makebox(0,0)[b]{$c$}}
      \put(45,33){\makebox(0,0)[b]{$d$}}
      \put(65,33){\makebox(0,0)[b]{$c$}}
      \put(85,33){\makebox(0,0)[b]{$b$}}
      \put(105,33){\makebox(0,0)[b]{$a$}}
    \end{picture}
    \begin{picture}(25,50)(0,5)
      \put(0,20){\vector(1,0){25}}
    \end{picture}
    \begin{picture}(120,50)(5,5)
      \put(25,20){\circle*{5}}
      \put(45,20){\circle*{5}}
      \put(65,20){\circle*{5}}
      \put(85,20){\circle*{5}}
      \put(105,20){\circle*{5}}
      \put(25,20){\line(1,0){80}}
      \qbezier(25,20)(65,-5)(105,20)
      \put(25,33){\makebox(0,0)[b]{$c$}}
      \put(45,33){\makebox(0,0)[b]{$d$}}
      \put(65,33){\makebox(0,0)[b]{$a$}}
      \put(85,33){\makebox(0,0)[b]{$b$}}
      \put(105,33){\makebox(0,0)[b]{$a$}}
    \end{picture}
    \begin{picture}(25,50)(0,5)
      \put(0,20){\vector(1,0){25}}
    \end{picture}
    \begin{picture}(120,50)(5,5)
      \put(25,20){\circle*{5}}
      \put(45,20){\circle*{5}}
      \put(65,20){\circle*{5}}
      \put(85,20){\circle*{5}}
      \put(105,20){\circle*{5}}
      \put(25,20){\line(1,0){80}}
      \qbezier(25,20)(65,-5)(105,20)
      \put(25,33){\makebox(0,0)[b]{$c$}}
      \put(45,33){\makebox(0,0)[b]{$d$}}
      \put(65,33){\makebox(0,0)[b]{$a$}}
      \put(85,33){\makebox(0,0)[b]{$b$}}
      \put(105,33){\makebox(0,0)[b]{$d$}}
    \end{picture}
    \begin{picture}(25,50)(0,5)
      \put(0,20){\vector(1,0){25}}
    \end{picture}
    \begin{picture}(120,50)(5,5)
      \put(25,20){\circle*{5}}
      \put(45,20){\circle*{5}}
      \put(65,20){\circle*{5}}
      \put(85,20){\circle*{5}}
      \put(105,20){\circle*{5}}
      \put(25,20){\line(1,0){80}}
      \qbezier(25,20)(65,-5)(105,20)
      \put(25,33){\makebox(0,0)[b]{$c$}}
      \put(45,33){\makebox(0,0)[b]{$b$}}
      \put(65,33){\makebox(0,0)[b]{$a$}}
      \put(85,33){\makebox(0,0)[b]{$b$}}
      \put(105,33){\makebox(0,0)[b]{$d$}}
    \end{picture}
    \begin{picture}(25,50)(0,5)
      \put(0,20){\vector(1,0){25}}
    \end{picture}
    \begin{picture}(120,50)(5,5)
      \put(25,20){\circle*{5}}
      \put(45,20){\circle*{5}}
      \put(65,20){\circle*{5}}
      \put(85,20){\circle*{5}}
      \put(105,20){\circle*{5}}
      \put(25,20){\line(1,0){80}}
      \qbezier(25,20)(65,-5)(105,20)
      \put(25,33){\makebox(0,0)[b]{$c$}}
      \put(45,33){\makebox(0,0)[b]{$b$}}
      \put(65,33){\makebox(0,0)[b]{$a$}}
      \put(85,33){\makebox(0,0)[b]{$c$}}
      \put(105,33){\makebox(0,0)[b]{$d$}}
    \end{picture}
    \begin{picture}(25,50)(0,5)
      \put(0,20){\vector(1,0){25}}
    \end{picture}
    \begin{picture}(120,50)(5,5)
      \put(25,20){\circle*{5}}
      \put(45,20){\circle*{5}}
      \put(65,20){\circle*{5}}
      \put(85,20){\circle*{5}}
      \put(105,20){\circle*{5}}
      \put(25,20){\line(1,0){80}}
      \qbezier(25,20)(65,-5)(105,20)
      \put(25,33){\makebox(0,0)[b]{$a$}}
      \put(45,33){\makebox(0,0)[b]{$b$}}
      \put(65,33){\makebox(0,0)[b]{$a$}}
      \put(85,33){\makebox(0,0)[b]{$c$}}
      \put(105,33){\makebox(0,0)[b]{$d$}}
    \end{picture}
  \end{center}

  By symmetry, we immediately get that $S_4(C_5)$ is a disjoint union of
  six copies of~$C_{20}$, so it is not connected.
\end{proof}

\begin{prop}\label{A.3}
  The strong colour graph $S_5(C_6)$ is connected.
\end{prop}

\begin{proof}
  In any strong 5-colouring of the 6-cycle~$C_6$, there are only two
  vertices having the same colour. Let~$\alpha$ be a strong 5-colouring of
  $C_6=v_1v_2...v_5v_6v_1$. We call~$\alpha$ an \emph{$a$-standard
    colouring} if $\alpha(v_1)=\alpha(v_3)=a$.

  \begin{figure}[ht]
    \centering
    \unitlength0.25mm
    \begin{picture}(145,50)(5,5)
      \put(25,20){\circle*{5}}
      \put(45,20){\circle*{5}}
      \put(65,20){\circle*{5}}
      \put(85,20){\circle*{5}}
      \put(105,20){\circle*{5}}
      \put(125,20){\circle*{5}}
      \put(25,20){\line(1,0){100}}
      \qbezier(25,20)(75,-5)(125,20)
      \put(25,33){\makebox(0,0)[b]{$a$}}
      \put(45,33){\makebox(0,0)[b]{$b$}}
      \put(65,33){\makebox(0,0)[b]{$a$}}
      \put(85,33){\makebox(0,0)[b]{$c$}}
      \put(105,33){\makebox(0,0)[b]{$d$}}
      \put(125,33){\makebox(0,0)[b]{$e$}}
      \put(25,10){\makebox(0,0)[t]{$v_1$}}
      \put(125,10){\makebox(0,0)[t]{$v_6$}}
    \end{picture}

    \caption{An $a$-standard colouring.}
    \label{f4}
  \end{figure}
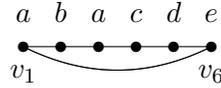

  We will prove the proposition by showing the following three steps.

  \medskip\noindent
  \textbf{Step 1}:\quad There is a path from $(a,b,a,c,d,e)$ to any other
  $a$-standard colourings.\\*[\parskip]
  First, we will show that there is a path from $(a,b,a,c,d,e)$ to
  $(a,c,a,b,d,e)$\,:

  \begin{center}
    \unitlength0.20mm
    \begin{picture}(140,50)(5,5)
      \put(25,20){\circle*{5}}
      \put(45,20){\circle*{5}}
      \put(65,20){\circle*{5}}
      \put(85,20){\circle*{5}}
      \put(105,20){\circle*{5}}
      \put(125,20){\circle*{5}}
      \put(25,20){\line(1,0){100}}
      \qbezier(25,20)(75,-5)(125,20)
      \put(25,33){\makebox(0,0)[b]{$a$}}
      \put(45,33){\makebox(0,0)[b]{$b$}}
      \put(65,33){\makebox(0,0)[b]{$a$}}
      \put(85,33){\makebox(0,0)[b]{$c$}}
      \put(105,33){\makebox(0,0)[b]{$d$}}
      \put(125,33){\makebox(0,0)[b]{$e$}}
    \end{picture}
    \begin{picture}(25,50)(0,5)
      \put(0,20){\vector(1,0){25}}
    \end{picture}
    \begin{picture}(140,50)(5,5)
      \put(25,20){\circle*{5}}
      \put(45,20){\circle*{5}}
      \put(65,20){\circle*{5}}
      \put(85,20){\circle*{5}}
      \put(105,20){\circle*{5}}
      \put(125,20){\circle*{5}}
      \put(25,20){\line(1,0){100}}
      \qbezier(25,20)(75,-5)(125,20)
      \put(25,33){\makebox(0,0)[b]{$a$}}
      \put(45,33){\makebox(0,0)[b]{$b$}}
      \put(65,33){\makebox(0,0)[b]{$e$}}
      \put(85,33){\makebox(0,0)[b]{$c$}}
      \put(105,33){\makebox(0,0)[b]{$d$}}
      \put(125,33){\makebox(0,0)[b]{$e$}}
    \end{picture}
    \begin{picture}(25,50)(0,5)
      \put(0,20){\vector(1,0){25}}
    \end{picture}
    \begin{picture}(140,50)(5,5)
      \put(25,20){\circle*{5}}
      \put(45,20){\circle*{5}}
      \put(65,20){\circle*{5}}
      \put(85,20){\circle*{5}}
      \put(105,20){\circle*{5}}
      \put(125,20){\circle*{5}}
      \put(25,20){\line(1,0){100}}
      \qbezier(25,20)(75,-5)(125,20)
      \put(25,33){\makebox(0,0)[b]{$a$}}
      \put(45,33){\makebox(0,0)[b]{$b$}}
      \put(65,33){\makebox(0,0)[b]{$e$}}
      \put(85,33){\makebox(0,0)[b]{$c$}}
      \put(105,33){\makebox(0,0)[b]{$d$}}
      \put(125,33){\makebox(0,0)[b]{$c$}}
    \end{picture}
    \begin{picture}(25,50)(0,5)
      \put(0,20){\vector(1,0){25}}
    \end{picture}
    \begin{picture}(140,50)(5,5)
      \put(25,20){\circle*{5}}
      \put(45,20){\circle*{5}}
      \put(65,20){\circle*{5}}
      \put(85,20){\circle*{5}}
      \put(105,20){\circle*{5}}
      \put(125,20){\circle*{5}}
      \put(25,20){\line(1,0){100}}
      \qbezier(25,20)(75,-5)(125,20)
      \put(25,33){\makebox(0,0)[b]{$a$}}
      \put(45,33){\makebox(0,0)[b]{$b$}}
      \put(65,33){\makebox(0,0)[b]{$e$}}
      \put(85,33){\makebox(0,0)[b]{$b$}}
      \put(105,33){\makebox(0,0)[b]{$d$}}
      \put(125,33){\makebox(0,0)[b]{$c$}}
    \end{picture}
    \begin{picture}(25,50)(0,5)
      \put(0,20){\vector(1,0){25}}
    \end{picture}
    \begin{picture}(140,50)(5,5)
      \put(25,20){\circle*{5}}
      \put(45,20){\circle*{5}}
      \put(65,20){\circle*{5}}
      \put(85,20){\circle*{5}}
      \put(105,20){\circle*{5}}
      \put(125,20){\circle*{5}}
      \put(25,20){\line(1,0){100}}
      \qbezier(25,20)(75,-5)(125,20)
      \put(25,33){\makebox(0,0)[b]{$a$}}
      \put(45,33){\makebox(0,0)[b]{$c$}}
      \put(65,33){\makebox(0,0)[b]{$e$}}
      \put(85,33){\makebox(0,0)[b]{$b$}}
      \put(105,33){\makebox(0,0)[b]{$d$}}
      \put(125,33){\makebox(0,0)[b]{$c$}}
    \end{picture}
    \begin{picture}(25,50)(0,5)
      \put(0,20){\vector(1,0){25}}
    \end{picture}
    \begin{picture}(140,50)(5,5)
      \put(25,20){\circle*{5}}
      \put(45,20){\circle*{5}}
      \put(65,20){\circle*{5}}
      \put(85,20){\circle*{5}}
      \put(105,20){\circle*{5}}
      \put(125,20){\circle*{5}}
      \put(25,20){\line(1,0){100}}
      \qbezier(25,20)(75,-5)(125,20)
      \put(25,33){\makebox(0,0)[b]{$a$}}
      \put(45,33){\makebox(0,0)[b]{$c$}}
      \put(65,33){\makebox(0,0)[b]{$e$}}
      \put(85,33){\makebox(0,0)[b]{$b$}}
      \put(105,33){\makebox(0,0)[b]{$d$}}
      \put(125,33){\makebox(0,0)[b]{$e$}}
    \end{picture}
    \begin{picture}(25,50)(0,5)
      \put(0,20){\vector(1,0){25}}
    \end{picture}
    \begin{picture}(140,50)(5,5)
      \put(25,20){\circle*{5}}
      \put(45,20){\circle*{5}}
      \put(65,20){\circle*{5}}
      \put(85,20){\circle*{5}}
      \put(105,20){\circle*{5}}
      \put(125,20){\circle*{5}}
      \put(25,20){\line(1,0){100}}
      \qbezier(25,20)(75,-5)(125,20)
      \put(25,33){\makebox(0,0)[b]{$a$}}
      \put(45,33){\makebox(0,0)[b]{$c$}}
      \put(65,33){\makebox(0,0)[b]{$a$}}
      \put(85,33){\makebox(0,0)[b]{$b$}}
      \put(105,33){\makebox(0,0)[b]{$d$}}
      \put(125,33){\makebox(0,0)[b]{$e$}}
    \end{picture}
  \end{center}

  By symmetry, there is also a path from $(a,b,a,c,d,e)$ to
  $(a,e,a,c,d,b)$.

  Next, we show that there is a path from $(a,b,a,c,d,e)$ to
  $(a,d,a,c,b,e)$\,:

  \begin{center}
    \unitlength0.20mm
    \begin{picture}(140,50)(5,5)
      \put(25,20){\circle*{5}}
      \put(45,20){\circle*{5}}
      \put(65,20){\circle*{5}}
      \put(85,20){\circle*{5}}
      \put(105,20){\circle*{5}}
      \put(125,20){\circle*{5}}
      \put(25,20){\line(1,0){100}}
      \qbezier(25,20)(75,-5)(125,20)
      \put(25,33){\makebox(0,0)[b]{$a$}}
      \put(45,33){\makebox(0,0)[b]{$b$}}
      \put(65,33){\makebox(0,0)[b]{$a$}}
      \put(85,33){\makebox(0,0)[b]{$c$}}
      \put(105,33){\makebox(0,0)[b]{$d$}}
      \put(125,33){\makebox(0,0)[b]{$e$}}
    \end{picture}
    \begin{picture}(25,50)(0,5)
      \put(0,20){\vector(1,0){25}}
    \end{picture}
    \begin{picture}(140,50)(5,5)
      \put(25,20){\circle*{5}}
      \put(45,20){\circle*{5}}
      \put(65,20){\circle*{5}}
      \put(85,20){\circle*{5}}
      \put(105,20){\circle*{5}}
      \put(125,20){\circle*{5}}
      \put(25,20){\line(1,0){100}}
      \qbezier(25,20)(75,-5)(125,20)
      \put(25,33){\makebox(0,0)[b]{$a$}}
      \put(45,33){\makebox(0,0)[b]{$b$}}
      \put(65,33){\makebox(0,0)[b]{$e$}}
      \put(85,33){\makebox(0,0)[b]{$c$}}
      \put(105,33){\makebox(0,0)[b]{$d$}}
      \put(125,33){\makebox(0,0)[b]{$e$}}
    \end{picture}
    \begin{picture}(25,50)(0,5)
      \put(0,20){\vector(1,0){25}}
    \end{picture}
    \begin{picture}(140,50)(5,5)
      \put(25,20){\circle*{5}}
      \put(45,20){\circle*{5}}
      \put(65,20){\circle*{5}}
      \put(85,20){\circle*{5}}
      \put(105,20){\circle*{5}}
      \put(125,20){\circle*{5}}
      \put(25,20){\line(1,0){100}}
      \qbezier(25,20)(75,-5)(125,20)
      \put(25,33){\makebox(0,0)[b]{$a$}}
      \put(45,33){\makebox(0,0)[b]{$b$}}
      \put(65,33){\makebox(0,0)[b]{$e$}}
      \put(85,33){\makebox(0,0)[b]{$c$}}
      \put(105,33){\makebox(0,0)[b]{$d$}}
      \put(125,33){\makebox(0,0)[b]{$b$}}
    \end{picture}
    \begin{picture}(25,50)(0,5)
      \put(0,20){\vector(1,0){25}}
    \end{picture}
    \begin{picture}(140,50)(5,5)
      \put(25,20){\circle*{5}}
      \put(45,20){\circle*{5}}
      \put(65,20){\circle*{5}}
      \put(85,20){\circle*{5}}
      \put(105,20){\circle*{5}}
      \put(125,20){\circle*{5}}
      \put(25,20){\line(1,0){100}}
      \qbezier(25,20)(75,-5)(125,20)
      \put(25,33){\makebox(0,0)[b]{$a$}}
      \put(45,33){\makebox(0,0)[b]{$d$}}
      \put(65,33){\makebox(0,0)[b]{$e$}}
      \put(85,33){\makebox(0,0)[b]{$c$}}
      \put(105,33){\makebox(0,0)[b]{$d$}}
      \put(125,33){\makebox(0,0)[b]{$b$}}
    \end{picture}
    \begin{picture}(25,50)(0,5)
      \put(0,20){\vector(1,0){25}}
    \end{picture}
    \begin{picture}(140,50)(5,5)
      \put(25,20){\circle*{5}}
      \put(45,20){\circle*{5}}
      \put(65,20){\circle*{5}}
      \put(85,20){\circle*{5}}
      \put(105,20){\circle*{5}}
      \put(125,20){\circle*{5}}
      \put(25,20){\line(1,0){100}}
      \qbezier(25,20)(75,-5)(125,20)
      \put(25,33){\makebox(0,0)[b]{$a$}}
      \put(45,33){\makebox(0,0)[b]{$d$}}
      \put(65,33){\makebox(0,0)[b]{$e$}}
      \put(85,33){\makebox(0,0)[b]{$c$}}
      \put(105,33){\makebox(0,0)[b]{$e$}}
      \put(125,33){\makebox(0,0)[b]{$b$}}
    \end{picture}
    \begin{picture}(25,50)(0,5)
      \put(0,20){\vector(1,0){25}}
    \end{picture}
    \begin{picture}(140,50)(5,5)
      \put(25,20){\circle*{5}}
      \put(45,20){\circle*{5}}
      \put(65,20){\circle*{5}}
      \put(85,20){\circle*{5}}
      \put(105,20){\circle*{5}}
      \put(125,20){\circle*{5}}
      \put(25,20){\line(1,0){100}}
      \qbezier(25,20)(75,-5)(125,20)
      \put(25,33){\makebox(0,0)[b]{$a$}}
      \put(45,33){\makebox(0,0)[b]{$d$}}
      \put(65,33){\makebox(0,0)[b]{$b$}}
      \put(85,33){\makebox(0,0)[b]{$c$}}
      \put(105,33){\makebox(0,0)[b]{$e$}}
      \put(125,33){\makebox(0,0)[b]{$b$}}
    \end{picture}
    \begin{picture}(25,50)(0,5)
      \put(0,20){\vector(1,0){25}}
    \end{picture}
    \begin{picture}(140,50)(5,5)
      \put(25,20){\circle*{5}}
      \put(45,20){\circle*{5}}
      \put(65,20){\circle*{5}}
      \put(85,20){\circle*{5}}
      \put(105,20){\circle*{5}}
      \put(125,20){\circle*{5}}
      \put(25,20){\line(1,0){100}}
      \qbezier(25,20)(75,-5)(125,20)
      \put(25,33){\makebox(0,0)[b]{$a$}}
      \put(45,33){\makebox(0,0)[b]{$d$}}
      \put(65,33){\makebox(0,0)[b]{$b$}}
      \put(85,33){\makebox(0,0)[b]{$c$}}
      \put(105,33){\makebox(0,0)[b]{$e$}}
      \put(125,33){\makebox(0,0)[b]{$c$}}
    \end{picture}
    \begin{picture}(25,50)(0,5)
      \put(0,20){\vector(1,0){25}}
    \end{picture}
    \begin{picture}(140,50)(5,5)
      \put(25,20){\circle*{5}}
      \put(45,20){\circle*{5}}
      \put(65,20){\circle*{5}}
      \put(85,20){\circle*{5}}
      \put(105,20){\circle*{5}}
      \put(125,20){\circle*{5}}
      \put(25,20){\line(1,0){100}}
      \qbezier(25,20)(75,-5)(125,20)
      \put(25,33){\makebox(0,0)[b]{$a$}}
      \put(45,33){\makebox(0,0)[b]{$d$}}
      \put(65,33){\makebox(0,0)[b]{$b$}}
      \put(85,33){\makebox(0,0)[b]{$a$}}
      \put(105,33){\makebox(0,0)[b]{$e$}}
      \put(125,33){\makebox(0,0)[b]{$c$}}
    \end{picture}
    \begin{picture}(25,50)(0,5)
      \put(0,20){\vector(1,0){25}}
    \end{picture}
    \begin{picture}(140,50)(5,5)
      \put(25,20){\circle*{5}}
      \put(45,20){\circle*{5}}
      \put(65,20){\circle*{5}}
      \put(85,20){\circle*{5}}
      \put(105,20){\circle*{5}}
      \put(125,20){\circle*{5}}
      \put(25,20){\line(1,0){100}}
      \qbezier(25,20)(75,-5)(125,20)
      \put(25,33){\makebox(0,0)[b]{$b$}}
      \put(45,33){\makebox(0,0)[b]{$d$}}
      \put(65,33){\makebox(0,0)[b]{$b$}}
      \put(85,33){\makebox(0,0)[b]{$a$}}
      \put(105,33){\makebox(0,0)[b]{$e$}}
      \put(125,33){\makebox(0,0)[b]{$c$}}
    \end{picture}
    \begin{picture}(25,50)(0,5)
      \put(0,20){\vector(1,0){25}}
    \end{picture}
    \begin{picture}(140,50)(5,5)
      \put(25,20){\circle*{5}}
      \put(45,20){\circle*{5}}
      \put(65,20){\circle*{5}}
      \put(85,20){\circle*{5}}
      \put(105,20){\circle*{5}}
      \put(125,20){\circle*{5}}
      \put(25,20){\line(1,0){100}}
      \qbezier(25,20)(75,-5)(125,20)
      \put(25,33){\makebox(0,0)[b]{$b$}}
      \put(45,33){\makebox(0,0)[b]{$d$}}
      \put(65,33){\makebox(0,0)[b]{$e$}}
      \put(85,33){\makebox(0,0)[b]{$a$}}
      \put(105,33){\makebox(0,0)[b]{$e$}}
      \put(125,33){\makebox(0,0)[b]{$c$}}
    \end{picture}
    \begin{picture}(25,50)(0,5)
      \put(0,20){\vector(1,0){25}}
    \end{picture}
    \begin{picture}(140,50)(5,5)
      \put(25,20){\circle*{5}}
      \put(45,20){\circle*{5}}
      \put(65,20){\circle*{5}}
      \put(85,20){\circle*{5}}
      \put(105,20){\circle*{5}}
      \put(125,20){\circle*{5}}
      \put(25,20){\line(1,0){100}}
      \qbezier(25,20)(75,-5)(125,20)
      \put(25,33){\makebox(0,0)[b]{$b$}}
      \put(45,33){\makebox(0,0)[b]{$d$}}
      \put(65,33){\makebox(0,0)[b]{$e$}}
      \put(85,33){\makebox(0,0)[b]{$a$}}
      \put(105,33){\makebox(0,0)[b]{$b$}}
      \put(125,33){\makebox(0,0)[b]{$c$}}
    \end{picture}
    \begin{picture}(25,50)(0,5)
      \put(0,20){\vector(1,0){25}}
    \end{picture}
    \begin{picture}(140,50)(5,5)
      \put(25,20){\circle*{5}}
      \put(45,20){\circle*{5}}
      \put(65,20){\circle*{5}}
      \put(85,20){\circle*{5}}
      \put(105,20){\circle*{5}}
      \put(125,20){\circle*{5}}
      \put(25,20){\line(1,0){100}}
      \qbezier(25,20)(75,-5)(125,20)
      \put(25,33){\makebox(0,0)[b]{$a$}}
      \put(45,33){\makebox(0,0)[b]{$d$}}
      \put(65,33){\makebox(0,0)[b]{$e$}}
      \put(85,33){\makebox(0,0)[b]{$a$}}
      \put(105,33){\makebox(0,0)[b]{$b$}}
      \put(125,33){\makebox(0,0)[b]{$c$}}
    \end{picture}
    \begin{picture}(25,50)(0,5)
      \put(0,20){\vector(1,0){25}}
    \end{picture}
    \begin{picture}(140,50)(5,5)
      \put(25,20){\circle*{5}}
      \put(45,20){\circle*{5}}
      \put(65,20){\circle*{5}}
      \put(85,20){\circle*{5}}
      \put(105,20){\circle*{5}}
      \put(125,20){\circle*{5}}
      \put(25,20){\line(1,0){100}}
      \qbezier(25,20)(75,-5)(125,20)
      \put(25,33){\makebox(0,0)[b]{$a$}}
      \put(45,33){\makebox(0,0)[b]{$d$}}
      \put(65,33){\makebox(0,0)[b]{$e$}}
      \put(85,33){\makebox(0,0)[b]{$c$}}
      \put(105,33){\makebox(0,0)[b]{$b$}}
      \put(125,33){\makebox(0,0)[b]{$c$}}
    \end{picture}
    \begin{picture}(25,50)(0,5)
      \put(0,20){\vector(1,0){25}}
    \end{picture}
    \begin{picture}(140,50)(5,5)
      \put(25,20){\circle*{5}}
      \put(45,20){\circle*{5}}
      \put(65,20){\circle*{5}}
      \put(85,20){\circle*{5}}
      \put(105,20){\circle*{5}}
      \put(125,20){\circle*{5}}
      \put(25,20){\line(1,0){100}}
      \qbezier(25,20)(75,-5)(125,20)
      \put(25,33){\makebox(0,0)[b]{$a$}}
      \put(45,33){\makebox(0,0)[b]{$d$}}
      \put(65,33){\makebox(0,0)[b]{$e$}}
      \put(85,33){\makebox(0,0)[b]{$c$}}
      \put(105,33){\makebox(0,0)[b]{$b$}}
      \put(125,33){\makebox(0,0)[b]{$e$}}
    \end{picture}
    \begin{picture}(25,50)(0,5)
      \put(0,20){\vector(1,0){25}}
    \end{picture}
    \begin{picture}(140,50)(5,5)
      \put(25,20){\circle*{5}}
      \put(45,20){\circle*{5}}
      \put(65,20){\circle*{5}}
      \put(85,20){\circle*{5}}
      \put(105,20){\circle*{5}}
      \put(125,20){\circle*{5}}
      \put(25,20){\line(1,0){100}}
      \qbezier(25,20)(75,-5)(125,20)
      \put(25,33){\makebox(0,0)[b]{$a$}}
      \put(45,33){\makebox(0,0)[b]{$d$}}
      \put(65,33){\makebox(0,0)[b]{$a$}}
      \put(85,33){\makebox(0,0)[b]{$c$}}
      \put(105,33){\makebox(0,0)[b]{$b$}}
      \put(125,33){\makebox(0,0)[b]{$e$}}
    \end{picture}
  \end{center}

  Now we consider $a$-standard colourings as permutations of $\{b,c,d,e\}$.
  Note that all these permutations can be generated by the transpositions
  $(b,c)$, $(b,e)$ and $(b,d)$. Therefore, since we can find a path from
  $(a,b,a,c,d,e)$ to $(a,c,a,b,d,e)$, from $(a,b,a,c,d,e)$ to
  $(a,e,a,c,d,b)$, and from $(a,b,a,c,d,e)$ to $(a,d,a,c,b,e)$, there is a
  path from $(e,a,e,b,c,d)$ to any other $a$-standard colourings.

  \medskip\noindent
  \textbf{Step 2}:\quad There is a path between any two types of standard
  colourings.\\*[\parskip]
  First, here is a path from an $a$-standard colouring~$\alpha$ to an
  $\alpha(v_5)$-standard colouring\,:

  \begin{center}
    \unitlength0.20mm
    \begin{picture}(140,50)(5,5)
      \put(25,20){\circle*{5}}
      \put(45,20){\circle*{5}}
      \put(65,20){\circle*{5}}
      \put(85,20){\circle*{5}}
      \put(105,20){\circle*{5}}
      \put(125,20){\circle*{5}}
      \put(25,20){\line(1,0){100}}
      \qbezier(25,20)(75,-5)(125,20)
      \put(25,33){\makebox(0,0)[b]{$a$}}
      \put(45,33){\makebox(0,0)[b]{$b$}}
      \put(65,33){\makebox(0,0)[b]{$a$}}
      \put(85,33){\makebox(0,0)[b]{$c$}}
      \put(105,33){\makebox(0,0)[b]{$d$}}
      \put(125,33){\makebox(0,0)[b]{$e$}}
    \end{picture}
    \begin{picture}(25,50)(0,5)
      \put(0,20){\vector(1,0){25}}
    \end{picture}
    \begin{picture}(140,50)(5,5)
      \put(25,20){\circle*{5}}
      \put(45,20){\circle*{5}}
      \put(65,20){\circle*{5}}
      \put(85,20){\circle*{5}}
      \put(105,20){\circle*{5}}
      \put(125,20){\circle*{5}}
      \put(25,20){\line(1,0){100}}
      \qbezier(25,20)(75,-5)(125,20)
      \put(25,33){\makebox(0,0)[b]{$a$}}
      \put(45,33){\makebox(0,0)[b]{$b$}}
      \put(65,33){\makebox(0,0)[b]{$d$}}
      \put(85,33){\makebox(0,0)[b]{$c$}}
      \put(105,33){\makebox(0,0)[b]{$d$}}
      \put(125,33){\makebox(0,0)[b]{$e$}}
    \end{picture}
    \begin{picture}(25,50)(0,5)
      \put(0,20){\vector(1,0){25}}
    \end{picture}
    \begin{picture}(140,50)(5,5)
      \put(25,20){\circle*{5}}
      \put(45,20){\circle*{5}}
      \put(65,20){\circle*{5}}
      \put(85,20){\circle*{5}}
      \put(105,20){\circle*{5}}
      \put(125,20){\circle*{5}}
      \put(25,20){\line(1,0){100}}
      \qbezier(25,20)(75,-5)(125,20)
      \put(25,33){\makebox(0,0)[b]{$a$}}
      \put(45,33){\makebox(0,0)[b]{$b$}}
      \put(65,33){\makebox(0,0)[b]{$d$}}
      \put(85,33){\makebox(0,0)[b]{$c$}}
      \put(105,33){\makebox(0,0)[b]{$a$}}
      \put(125,33){\makebox(0,0)[b]{$e$}}
    \end{picture}
    \begin{picture}(25,50)(0,5)
      \put(0,20){\vector(1,0){25}}
    \end{picture}
    \begin{picture}(140,50)(5,5)
      \put(25,20){\circle*{5}}
      \put(45,20){\circle*{5}}
      \put(65,20){\circle*{5}}
      \put(85,20){\circle*{5}}
      \put(105,20){\circle*{5}}
      \put(125,20){\circle*{5}}
      \put(25,20){\line(1,0){100}}
      \qbezier(25,20)(75,-5)(125,20)
      \put(25,33){\makebox(0,0)[b]{$d$}}
      \put(45,33){\makebox(0,0)[b]{$b$}}
      \put(65,33){\makebox(0,0)[b]{$d$}}
      \put(85,33){\makebox(0,0)[b]{$c$}}
      \put(105,33){\makebox(0,0)[b]{$a$}}
      \put(125,33){\makebox(0,0)[b]{$e$}}
    \end{picture}
  \end{center}

  Next, a path from an $a$-standard colouring~$\alpha$ to an
  $\alpha(v_4)$-standard colouring\,:

  \begin{center}
    \unitlength0.20mm
    \begin{picture}(140,50)(5,5)
      \put(25,20){\circle*{5}}
      \put(45,20){\circle*{5}}
      \put(65,20){\circle*{5}}
      \put(85,20){\circle*{5}}
      \put(105,20){\circle*{5}}
      \put(125,20){\circle*{5}}
      \put(25,20){\line(1,0){100}}
      \qbezier(25,20)(75,-5)(125,20)
      \put(25,33){\makebox(0,0)[b]{$a$}}
      \put(45,33){\makebox(0,0)[b]{$b$}}
      \put(65,33){\makebox(0,0)[b]{$a$}}
      \put(85,33){\makebox(0,0)[b]{$c$}}
      \put(105,33){\makebox(0,0)[b]{$d$}}
      \put(125,33){\makebox(0,0)[b]{$e$}}
    \end{picture}
    \begin{picture}(25,50)(0,5)
      \put(0,20){\vector(1,0){25}}
    \end{picture}
    \begin{picture}(140,50)(5,5)
      \put(25,20){\circle*{5}}
      \put(45,20){\circle*{5}}
      \put(65,20){\circle*{5}}
      \put(85,20){\circle*{5}}
      \put(105,20){\circle*{5}}
      \put(125,20){\circle*{5}}
      \put(25,20){\line(1,0){100}}
      \qbezier(25,20)(75,-5)(125,20)
      \put(25,33){\makebox(0,0)[b]{$c$}}
      \put(45,33){\makebox(0,0)[b]{$b$}}
      \put(65,33){\makebox(0,0)[b]{$a$}}
      \put(85,33){\makebox(0,0)[b]{$c$}}
      \put(105,33){\makebox(0,0)[b]{$d$}}
      \put(125,33){\makebox(0,0)[b]{$e$}}
    \end{picture}
    \begin{picture}(25,50)(0,5)
      \put(0,20){\vector(1,0){25}}
    \end{picture}
    \begin{picture}(140,50)(5,5)
      \put(25,20){\circle*{5}}
      \put(45,20){\circle*{5}}
      \put(65,20){\circle*{5}}
      \put(85,20){\circle*{5}}
      \put(105,20){\circle*{5}}
      \put(125,20){\circle*{5}}
      \put(25,20){\line(1,0){100}}
      \qbezier(25,20)(75,-5)(125,20)
      \put(25,33){\makebox(0,0)[b]{$c$}}
      \put(45,33){\makebox(0,0)[b]{$b$}}
      \put(65,33){\makebox(0,0)[b]{$a$}}
      \put(85,33){\makebox(0,0)[b]{$e$}}
      \put(105,33){\makebox(0,0)[b]{$d$}}
      \put(125,33){\makebox(0,0)[b]{$e$}}
    \end{picture}
    \begin{picture}(25,50)(0,5)
      \put(0,20){\vector(1,0){25}}
    \end{picture}
    \begin{picture}(140,50)(5,5)
      \put(25,20){\circle*{5}}
      \put(45,20){\circle*{5}}
      \put(65,20){\circle*{5}}
      \put(85,20){\circle*{5}}
      \put(105,20){\circle*{5}}
      \put(125,20){\circle*{5}}
      \put(25,20){\line(1,0){100}}
      \qbezier(25,20)(75,-5)(125,20)
      \put(25,33){\makebox(0,0)[b]{$c$}}
      \put(45,33){\makebox(0,0)[b]{$b$}}
      \put(65,33){\makebox(0,0)[b]{$a$}}
      \put(85,33){\makebox(0,0)[b]{$e$}}
      \put(105,33){\makebox(0,0)[b]{$d$}}
      \put(125,33){\makebox(0,0)[b]{$a$}}
    \end{picture}
    \begin{picture}(25,50)(0,5)
      \put(0,20){\vector(1,0){25}}
    \end{picture}
    \begin{picture}(140,50)(5,5)
      \put(25,20){\circle*{5}}
      \put(45,20){\circle*{5}}
      \put(65,20){\circle*{5}}
      \put(85,20){\circle*{5}}
      \put(105,20){\circle*{5}}
      \put(125,20){\circle*{5}}
      \put(25,20){\line(1,0){100}}
      \qbezier(25,20)(75,-5)(125,20)
      \put(25,33){\makebox(0,0)[b]{$c$}}
      \put(45,33){\makebox(0,0)[b]{$b$}}
      \put(65,33){\makebox(0,0)[b]{$c$}}
      \put(85,33){\makebox(0,0)[b]{$e$}}
      \put(105,33){\makebox(0,0)[b]{$d$}}
      \put(125,33){\makebox(0,0)[b]{$a$}}
    \end{picture}
  \end{center}

  By symmetry, there is also a path from an $a$-standard colouring~$\alpha$
  to an $\alpha(v_6)$-standard colouring.

  And finally, a path from an $a$-standard colouring~$\alpha$ to an
  $\alpha(v_2)$-standard colouring\,:

  \begin{center}
    \unitlength0.20mm
    \begin{picture}(140,50)(5,5)
      \put(25,20){\circle*{5}}
      \put(45,20){\circle*{5}}
      \put(65,20){\circle*{5}}
      \put(85,20){\circle*{5}}
      \put(105,20){\circle*{5}}
      \put(125,20){\circle*{5}}
      \put(25,20){\line(1,0){100}}
      \qbezier(25,20)(75,-5)(125,20)
      \put(25,33){\makebox(0,0)[b]{$a$}}
      \put(45,33){\makebox(0,0)[b]{$b$}}
      \put(65,33){\makebox(0,0)[b]{$a$}}
      \put(85,33){\makebox(0,0)[b]{$c$}}
      \put(105,33){\makebox(0,0)[b]{$d$}}
      \put(125,33){\makebox(0,0)[b]{$e$}}
    \end{picture}
    \begin{picture}(25,50)(0,5)
      \put(0,20){\vector(1,0){25}}
    \end{picture}
    \begin{picture}(140,50)(5,5)
      \put(25,20){\circle*{5}}
      \put(45,20){\circle*{5}}
      \put(65,20){\circle*{5}}
      \put(85,20){\circle*{5}}
      \put(105,20){\circle*{5}}
      \put(125,20){\circle*{5}}
      \put(25,20){\line(1,0){100}}
      \qbezier(25,20)(75,-5)(125,20)
      \put(25,33){\makebox(0,0)[b]{$d$}}
      \put(45,33){\makebox(0,0)[b]{$b$}}
      \put(65,33){\makebox(0,0)[b]{$a$}}
      \put(85,33){\makebox(0,0)[b]{$c$}}
      \put(105,33){\makebox(0,0)[b]{$d$}}
      \put(125,33){\makebox(0,0)[b]{$e$}}
    \end{picture}
    \begin{picture}(25,50)(0,5)
      \put(0,20){\vector(1,0){25}}
    \end{picture}
    \begin{picture}(140,50)(5,5)
      \put(25,20){\circle*{5}}
      \put(45,20){\circle*{5}}
      \put(65,20){\circle*{5}}
      \put(85,20){\circle*{5}}
      \put(105,20){\circle*{5}}
      \put(125,20){\circle*{5}}
      \put(25,20){\line(1,0){100}}
      \qbezier(25,20)(75,-5)(125,20)
      \put(25,33){\makebox(0,0)[b]{$d$}}
      \put(45,33){\makebox(0,0)[b]{$b$}}
      \put(65,33){\makebox(0,0)[b]{$a$}}
      \put(85,33){\makebox(0,0)[b]{$c$}}
      \put(105,33){\makebox(0,0)[b]{$b$}}
      \put(125,33){\makebox(0,0)[b]{$e$}}
    \end{picture}
    \begin{picture}(25,50)(0,5)
      \put(0,20){\vector(1,0){25}}
    \end{picture}
    \begin{picture}(140,50)(5,5)
      \put(25,20){\circle*{5}}
      \put(45,20){\circle*{5}}
      \put(65,20){\circle*{5}}
      \put(85,20){\circle*{5}}
      \put(105,20){\circle*{5}}
      \put(125,20){\circle*{5}}
      \put(25,20){\line(1,0){100}}
      \qbezier(25,20)(75,-5)(125,20)
      \put(25,33){\makebox(0,0)[b]{$d$}}
      \put(45,33){\makebox(0,0)[b]{$c$}}
      \put(65,33){\makebox(0,0)[b]{$a$}}
      \put(85,33){\makebox(0,0)[b]{$c$}}
      \put(105,33){\makebox(0,0)[b]{$b$}}
      \put(125,33){\makebox(0,0)[b]{$e$}}
    \end{picture}
    \begin{picture}(25,50)(0,5)
      \put(0,20){\vector(1,0){25}}
    \end{picture}
    \begin{picture}(140,50)(5,5)
      \put(25,20){\circle*{5}}
      \put(45,20){\circle*{5}}
      \put(65,20){\circle*{5}}
      \put(85,20){\circle*{5}}
      \put(105,20){\circle*{5}}
      \put(125,20){\circle*{5}}
      \put(25,20){\line(1,0){100}}
      \qbezier(25,20)(75,-5)(125,20)
      \put(25,33){\makebox(0,0)[b]{$d$}}
      \put(45,33){\makebox(0,0)[b]{$c$}}
      \put(65,33){\makebox(0,0)[b]{$a$}}
      \put(85,33){\makebox(0,0)[b]{$d$}}
      \put(105,33){\makebox(0,0)[b]{$b$}}
      \put(125,33){\makebox(0,0)[b]{$e$}}
    \end{picture}
    \begin{picture}(25,50)(0,5)
      \put(0,20){\vector(1,0){25}}
    \end{picture}
    \begin{picture}(140,50)(5,5)
      \put(25,20){\circle*{5}}
      \put(45,20){\circle*{5}}
      \put(65,20){\circle*{5}}
      \put(85,20){\circle*{5}}
      \put(105,20){\circle*{5}}
      \put(125,20){\circle*{5}}
      \put(25,20){\line(1,0){100}}
      \qbezier(25,20)(75,-5)(125,20)
      \put(25,33){\makebox(0,0)[b]{$b$}}
      \put(45,33){\makebox(0,0)[b]{$c$}}
      \put(65,33){\makebox(0,0)[b]{$a$}}
      \put(85,33){\makebox(0,0)[b]{$d$}}
      \put(105,33){\makebox(0,0)[b]{$b$}}
      \put(125,33){\makebox(0,0)[b]{$e$}}
    \end{picture}
    \begin{picture}(25,50)(0,5)
      \put(0,20){\vector(1,0){25}}
    \end{picture}
    \begin{picture}(140,50)(5,5)
      \put(25,20){\circle*{5}}
      \put(45,20){\circle*{5}}
      \put(65,20){\circle*{5}}
      \put(85,20){\circle*{5}}
      \put(105,20){\circle*{5}}
      \put(125,20){\circle*{5}}
      \put(25,20){\line(1,0){100}}
      \qbezier(25,20)(75,-5)(125,20)
      \put(25,33){\makebox(0,0)[b]{$b$}}
      \put(45,33){\makebox(0,0)[b]{$c$}}
      \put(65,33){\makebox(0,0)[b]{$a$}}
      \put(85,33){\makebox(0,0)[b]{$d$}}
      \put(105,33){\makebox(0,0)[b]{$a$}}
      \put(125,33){\makebox(0,0)[b]{$e$}}
    \end{picture}
    \begin{picture}(25,50)(0,5)
      \put(0,20){\vector(1,0){25}}
    \end{picture}
    \begin{picture}(140,50)(5,5)
      \put(25,20){\circle*{5}}
      \put(45,20){\circle*{5}}
      \put(65,20){\circle*{5}}
      \put(85,20){\circle*{5}}
      \put(105,20){\circle*{5}}
      \put(125,20){\circle*{5}}
      \put(25,20){\line(1,0){100}}
      \qbezier(25,20)(75,-5)(125,20)
      \put(25,33){\makebox(0,0)[b]{$b$}}
      \put(45,33){\makebox(0,0)[b]{$c$}}
      \put(65,33){\makebox(0,0)[b]{$b$}}
      \put(85,33){\makebox(0,0)[b]{$d$}}
      \put(105,33){\makebox(0,0)[b]{$a$}}
      \put(125,33){\makebox(0,0)[b]{$e$}}
    \end{picture}
  \end{center}

  \medskip\noindent
  \textbf{Step 3}:\quad Each colouring has a path to some standard
  colouring.\\*[\parskip]
  Let~$\alpha$ be a strong 5-colouring of~$C_6$. Then~$\alpha$ has one of
  the following forms\,:

  \begin{center}
    \unitlength0.20mm
    \begin{picture}(140,50)(5,5)
      \put(25,20){\circle*{5}}
      \put(45,20){\circle*{5}}
      \put(65,20){\circle*{5}}
      \put(85,20){\circle*{5}}
      \put(105,20){\circle*{5}}
      \put(125,20){\circle*{5}}
      \put(25,20){\line(1,0){100}}
      \qbezier(25,20)(75,-5)(125,20)
      \put(25,33){\makebox(0,0)[b]{$a$}}
      \put(45,33){\makebox(0,0)[b]{$b$}}
      \put(65,33){\makebox(0,0)[b]{$a$}}
      \put(85,33){\makebox(0,0)[b]{$c$}}
      \put(105,33){\makebox(0,0)[b]{$d$}}
      \put(125,33){\makebox(0,0)[b]{$e$}}
    \end{picture}
    \begin{picture}(140,50)(5,5)
      \put(25,20){\circle*{5}}
      \put(45,20){\circle*{5}}
      \put(65,20){\circle*{5}}
      \put(85,20){\circle*{5}}
      \put(105,20){\circle*{5}}
      \put(125,20){\circle*{5}}
      \put(25,20){\line(1,0){100}}
      \qbezier(25,20)(75,-5)(125,20)
      \put(25,33){\makebox(0,0)[b]{$a$}}
      \put(45,33){\makebox(0,0)[b]{$b$}}
      \put(65,33){\makebox(0,0)[b]{$c$}}
      \put(85,33){\makebox(0,0)[b]{$a$}}
      \put(105,33){\makebox(0,0)[b]{$d$}}
      \put(125,33){\makebox(0,0)[b]{$e$}}
    \end{picture}
    \begin{picture}(140,50)(5,5)
      \put(25,20){\circle*{5}}
      \put(45,20){\circle*{5}}
      \put(65,20){\circle*{5}}
      \put(85,20){\circle*{5}}
      \put(105,20){\circle*{5}}
      \put(125,20){\circle*{5}}
      \put(25,20){\line(1,0){100}}
      \qbezier(25,20)(75,-5)(125,20)
      \put(25,33){\makebox(0,0)[b]{$a$}}
      \put(45,33){\makebox(0,0)[b]{$b$}}
      \put(65,33){\makebox(0,0)[b]{$c$}}
      \put(85,33){\makebox(0,0)[b]{$d$}}
      \put(105,33){\makebox(0,0)[b]{$a$}}
      \put(125,33){\makebox(0,0)[b]{$e$}}
    \end{picture}
    \begin{picture}(140,50)(5,5)
      \put(25,20){\circle*{5}}
      \put(45,20){\circle*{5}}
      \put(65,20){\circle*{5}}
      \put(85,20){\circle*{5}}
      \put(105,20){\circle*{5}}
      \put(125,20){\circle*{5}}
      \put(25,20){\line(1,0){100}}
      \qbezier(25,20)(75,-5)(125,20)
      \put(25,33){\makebox(0,0)[b]{$b$}}
      \put(45,33){\makebox(0,0)[b]{$a$}}
      \put(65,33){\makebox(0,0)[b]{$c$}}
      \put(85,33){\makebox(0,0)[b]{$a$}}
      \put(105,33){\makebox(0,0)[b]{$d$}}
      \put(125,33){\makebox(0,0)[b]{$e$}}
    \end{picture}
    \begin{picture}(140,50)(5,5)
      \put(25,20){\circle*{5}}
      \put(45,20){\circle*{5}}
      \put(65,20){\circle*{5}}
      \put(85,20){\circle*{5}}
      \put(105,20){\circle*{5}}
      \put(125,20){\circle*{5}}
      \put(25,20){\line(1,0){100}}
      \qbezier(25,20)(75,-5)(125,20)
      \put(25,33){\makebox(0,0)[b]{$b$}}
      \put(45,33){\makebox(0,0)[b]{$a$}}
      \put(65,33){\makebox(0,0)[b]{$c$}}
      \put(85,33){\makebox(0,0)[b]{$d$}}
      \put(105,33){\makebox(0,0)[b]{$a$}}
      \put(125,33){\makebox(0,0)[b]{$e$}}
    \end{picture}
    \begin{picture}(140,50)(5,5)
      \put(25,20){\circle*{5}}
      \put(45,20){\circle*{5}}
      \put(65,20){\circle*{5}}
      \put(85,20){\circle*{5}}
      \put(105,20){\circle*{5}}
      \put(125,20){\circle*{5}}
      \put(25,20){\line(1,0){100}}
      \qbezier(25,20)(75,-5)(125,20)
      \put(25,33){\makebox(0,0)[b]{$b$}}
      \put(45,33){\makebox(0,0)[b]{$a$}}
      \put(65,33){\makebox(0,0)[b]{$c$}}
      \put(85,33){\makebox(0,0)[b]{$d$}}
      \put(105,33){\makebox(0,0)[b]{$e$}}
      \put(125,33){\makebox(0,0)[b]{$a$}}
    \end{picture}
    \begin{picture}(140,50)(5,5)
      \put(25,20){\circle*{5}}
      \put(45,20){\circle*{5}}
      \put(65,20){\circle*{5}}
      \put(85,20){\circle*{5}}
      \put(105,20){\circle*{5}}
      \put(125,20){\circle*{5}}
      \put(25,20){\line(1,0){100}}
      \qbezier(25,20)(75,-5)(125,20)
      \put(25,33){\makebox(0,0)[b]{$b$}}
      \put(45,33){\makebox(0,0)[b]{$c$}}
      \put(65,33){\makebox(0,0)[b]{$a$}}
      \put(85,33){\makebox(0,0)[b]{$d$}}
      \put(105,33){\makebox(0,0)[b]{$a$}}
      \put(125,33){\makebox(0,0)[b]{$e$}}
    \end{picture}
    \begin{picture}(140,50)(5,5)
      \put(25,20){\circle*{5}}
      \put(45,20){\circle*{5}}
      \put(65,20){\circle*{5}}
      \put(85,20){\circle*{5}}
      \put(105,20){\circle*{5}}
      \put(125,20){\circle*{5}}
      \put(25,20){\line(1,0){100}}
      \qbezier(25,20)(75,-5)(125,20)
      \put(25,33){\makebox(0,0)[b]{$b$}}
      \put(45,33){\makebox(0,0)[b]{$c$}}
      \put(65,33){\makebox(0,0)[b]{$a$}}
      \put(85,33){\makebox(0,0)[b]{$d$}}
      \put(105,33){\makebox(0,0)[b]{$e$}}
      \put(125,33){\makebox(0,0)[b]{$a$}}
    \end{picture}
    \begin{picture}(140,50)(5,5)
      \put(25,20){\circle*{5}}
      \put(45,20){\circle*{5}}
      \put(65,20){\circle*{5}}
      \put(85,20){\circle*{5}}
      \put(105,20){\circle*{5}}
      \put(125,20){\circle*{5}}
      \put(25,20){\line(1,0){100}}
      \qbezier(25,20)(75,-5)(125,20)
      \put(25,33){\makebox(0,0)[b]{$b$}}
      \put(45,33){\makebox(0,0)[b]{$c$}}
      \put(65,33){\makebox(0,0)[b]{$d$}}
      \put(85,33){\makebox(0,0)[b]{$a$}}
      \put(105,33){\makebox(0,0)[b]{$e$}}
      \put(125,33){\makebox(0,0)[b]{$a$}}
    \end{picture}
  \end{center}

  The first form is already an $a$-standard colouring. For the second and
  the third forms, we just recolour vertex~$v_1$ to~$c$, and for the
  seventh and eight forms, we just recolour vertex~$v_3$ to~$b$. For all
  the remaining colourings, we can find a path of length two to some
  standard colouring. We will leave checking that to the reader.

  \medskip\noindent
  It is straightforward to see that appropriate renaming of the colours and
  sequence of the paths in Steps 1\,--\,3 will transform any strong
  5-colouring of~$P_6$ into any other strong 5-colouring.
\end{proof}

\begin{thm}\label{A.4}
  The strong colour graph $S_k(C_n)$ is connected if and only if $k\ge4$,
  $n\ge6$, and $n\ge k+1$.
\end{thm}

\begin{proof}
  We already have seen that $S_3(C_n)$, $n\ge3$ and $S_4(C_5)$ are
  disconnected. From Theorems~\ref{Z.6} and~\ref{B.10} we easily obtain
  that $S_k(C_n)$ is connected for all $k\ge4$, $n\ge6$, and $n\ge k+2$.
  Since $S_5(C_6)$ is connected, all that is left to prove is that
  $S_k(C_{k+1})$ is connected for all $k\ge6$.

  Let $k\ge6$ and let~$\alpha$ and~$\beta$ be strong $k$-colourings of
  $C_{k+1}=v_1v_2\ldots v_{k+1}v_1$. In~$\alpha$, there will be a vertex,
  say~$v_1$, which has an unique colour, say colour~$a$. We say that a
  strong $k$-colouring of $C_{k+1}$ is \emph{good} if~$v_1$ is the only
  vertex in $C_{k+1}$, which is coloured~$a$. Thus~$\alpha$ is good.

  If~$\beta$ is good as well, then remove~$v_1$, and let~$\alpha'$
  and~$\beta'$ be the strong $(k-1)$-colourings of~$P_k$ obtained
  from~$\alpha$ and~$\beta$, respectively. Since $k\ge6$, by
  Proposition~\ref{B.9} there is a path in $S_{k-1}(P_k)$ from~$\alpha'$
  to~$\beta'$. Using the same steps gives a path from~$\alpha$ to~$\beta$
  in~$C_{k+1}$.

  So suppose that~$\beta$ is not good. As we colour the $k+1$ vertices of
  $C_{k+1}$ with~$k$ colours, there are only two vertices having the same
  colour. We distinguish five cases.

  \medskip\noindent
  \textbf{Case 1}:\quad In~$\beta$, $v_1$ is coloured~$a$, but there is
  a second vertex~$v_i$ coloured~$a$ as well.\\*[\parskip]
  Then just recolour~$v_i$ to another colour. The resulting colouring is
  good, and we are done by the paragraph above.

  \medskip\noindent
  \textbf{Case 2}:\quad In~$\beta$, $v_1$ and some other vertex~$v_i$ have
  the same colour $b\ne a$, while a third vertex~$v_j$ with $j\ne 2,k+1$ is
  coloured~$a$.\\*[\parskip]
  Then we first recolour~$v_1$ to~$a$, and then recolour~$v_j$ to another
  colour. Again, this gives a good colouring, so we are done.

  \medskip\noindent
  \textbf{Case 3}:\quad In~$\beta$, $v_1$ and some other vertex~$v_i$ have
  the same colour $b\ne a$, while a third vertex~$v_j$ with $j\in\{2,k+1\}$
  is coloured~$a$.\\*[\parskip]
  Without loss of generality, assume that $j=k+1$.

  \medskip\noindent
  \textbf{Subcase 3.1}:\quad We have $i\ge5$.\\*[\parskip]
  Now first recolour vertex~$v_i$ to colour~$\beta(v_3)$, then
  recolour~$v_3$~to $\beta(v_{k+1})=a$, $v_{k+1}$ to~$\beta(v_4)$, $v_4$
  to~$\beta(v_1)=b$, $v_1$ to~$a$, and finally recolour~$v_3$ to some
  colour different from~$a$. It is easy to check that the remaining
  colouring is good.

  \medskip\noindent
  \textbf{Subcase 3.2}:\quad We have $i\in\{3,4\}$.\\*[\parskip]
  Recall that $j=k+1\ge7$. Now first recolour vertex~$v_i$ to colour
  $\beta(v_6)$, then recolour~$v_6$~to $\beta(v_1)=b$. Now~$v_1$ and~$v_6$
  have the same colour~$b$, so we are back to Subcase~3.1.

  \medskip\noindent
  \textbf{Case 4}:\quad In~$\beta$, $v_1$ has a unique colour $b\ne a$,
  while there are two vertices~$v_i$ and~$v_j$ coloured~$a$.

  \medskip\noindent
  \textbf{Subcase 4.1}:\quad We have $\{i,j\}=\{2,k+1\}$.\\*[\parskip]
  Now first recolour vertex~$v_2$ to~$\beta(v_4)$, and then recolour~$v_4$
  to $\beta(v_1)=b$. Then we are back to Subcase~3.2.

  \medskip\noindent
  \textbf{Subcase 4.2}:\quad We have $\{i,j\}\ne\{2,k+1\}$.\\*[\parskip]
  Without loss of generality, assume that $i\ne2,k+1$. Then we can
  recolour~$v_i$ to $\beta(v_1)=b$. This means that~$v_1$ and~$v_i$ have
  the same colour $b\ne a$, so we are back in Case~2 or~3.

  \medskip\noindent
  \textbf{Case 5}:\quad In~$\beta$, $v_1$ has a unique colour $b\ne a$,
  there is a unique vertex~$v_i$ coloured~$a$, and two vertices~$v_j$
  and~$v_\ell$ have the same colour $c\ne a,b$.

  \medskip\noindent
  \textbf{Subcase 5.1}:\quad We have $\{j,\ell\}=\{2,k+1\}$.\\*[\parskip]
  Since $k+1\ge7$, we must have $i\ne3$ or $i\ne k$. Without loss of
  generality, assume that $i\ne3$. Then recolour~$v_2$ to $\beta(v_i)=a$,
  and next recolour~$v_i$ to $\beta(v_1)=b$. This brings us back to Case~3.

  \medskip\noindent
  \textbf{Subcase 5.2}:\quad We have $\{j,\ell\}\ne\{2,k+1\}$.\\*[\parskip]
  Without loss of generality, assume that $j\ne2,k+1$. Then we can
  recolour~$v_j$ to $\beta(v_1)=b$. This means that~$v_1$ and~$v_j$ have
  the same colour $b\ne a$, and we are back in Case~2 or~3.
\end{proof}

\section{The Strong 3-Colour Graph of Trees}\label{sec5}

  The aim of this section is to classify the trees~$T$ for which the
  strong 3-colour graph $S_3(T)$ is connected.

  For this we need to consider some special trees. First, in
  Section~\ref{sec2} we saw that the strong $k$-colour graph of a complete
  bipartite graph is not connected, so $S_3(K_{1,n})$ is disconnected for all
  $n\ge2$.

\begin{figure}[ht]
  \centering
  \unitlength0.28mm
  \begin{picture}(60,40)(15,-60)
    \put(15,0){\circle*{5}}\put(15,40){\circle*{5}}
    \put(45,0){\circle*{5}}\put(45,40){\circle*{5}}
    \put(75,0){\circle*{5}}\put(75,40){\circle*{5}}
    \put(15,0){\line(1,0){60}}
    \put(15,40){\line(1,0){60}}
    \put(45,0){\line(0,1){40}}
    \put(15,50){\makebox(0,0)[b]{$x_1$}}
    \put(45,50){\makebox(0,0)[b]{$x_2$}}
    \put(75,50){\makebox(0,0)[b]{$x_3$}}
    \put(15,-10){\makebox(0,0)[t]{$x_4$}}
    \put(45,-10){\makebox(0,0)[t]{$x_5$}}
    \put(75,-10){\makebox(0,0)[t]{$x_6$}}
  \end{picture}\qquad\qquad\qquad
  \begin{picture}(70,140)(25,-30)
    \put(25,70){\circle*{5}}
    \put(45,70){\circle*{5}}
    \put(65,70){\circle*{1}}
    \put(70,70){\circle*{1}}
    \put(75,70){\circle*{1}}
    \put(95,70){\circle*{5}}
    \put(65,40){\circle*{5}}
    \put(65,20){\circle*{5}}
    \put(65,0){\circle*{5}}
    \put(65,0){\line(0,1){40}}
    \put(65,40){\line(-4,3){40}}
    \put(65,40){\line(-2,3){20}}
    \put(65,40){\line(1,1){30}}
    \put(80,35){\makebox(0,0){$v_0$}}
    \put(25,80){\makebox(0,0)[b]{$v_1$}}
    \put(45,80){\makebox(0,0)[b]{$v_2$}}
    \put(95,80){\makebox(0,0)[b]{$v_n$}}
  \end{picture}\qquad\qquad\qquad
  \begin{picture}(70,140)(25,-30)
    \put(25,70){\circle*{5}}
    \put(45,70){\circle*{5}}
    \put(65,70){\circle*{1}}
    \put(70,70){\circle*{1}}
    \put(75,70){\circle*{1}}
    \put(95,70){\circle*{5}}
    \put(65,40){\circle*{5}}
    \put(65,20){\circle*{5}}
    \put(25,-10){\circle*{5}}
    \put(45,-10){\circle*{5}}
    \put(65,-10){\circle*{1}}
    \put(70,-10){\circle*{1}}
    \put(75,-10){\circle*{1}}
    \put(95,-10){\circle*{5}}
    \put(65,20){\line(0,1){20}}
    \put(65,40){\line(-4,3){40}}
    \put(65,40){\line(-2,3){20}}
    \put(65,40){\line(1,1){30}}
    \put(65,20){\line(-4,-3){40}}
    \put(65,20){\line(-2,-3){20}}
    \put(65,20){\line(1,-1){30}}
    \put(25,80){\makebox(0,0)[b]{$u_1$}}
    \put(45,80){\makebox(0,0)[b]{$u_2$}}
    \put(95,78){\makebox(0,0)[b]{$u_p$}}
    \put(25,-20){\makebox(0,0)[t]{$w_1$}}
    \put(45,-20){\makebox(0,0)[t]{$w_2$}}
    \put(95,-20){\makebox(0,0)[t]{$w_q$}}
  \end{picture}

  \caption{The graphs $I$, $\Psi_n$, and $\Phi_{p,q}$.}
  \label{f5}
\end{figure}
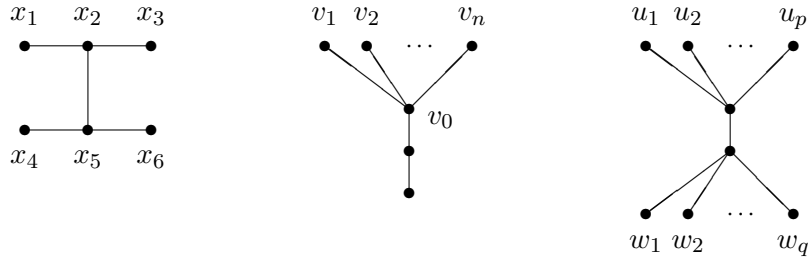

  For $n\ge1$ and $p,q\ge2$, let~$I$, $\Psi_n$ and $\Phi_{p,q}$ be the graphs
  sketched in Figure~\ref{f5}, respectively.

  It is straightforward to check that in any strong 3-colouring of~$\Psi_n$
  we cannot recolour the vertex~$v_0$ to another colour so that the
  resulting 3-colouring is strong again. Hence the strong colour graph
  $S_3(\Psi_n)$ is disconnected for all $n\ge1$.

\begin{prop} \label{D.5}
  The strong 3-colour graph $S_3(I)$ is connected.
\end{prop}

\begin{proof}
  Let~$\alpha$ be a strong 3-colouring of the graph~$I$, with vertex set
  $\{x_1,x_2,\ldots,x_6\}$ as in Figure~\ref{f5}. We call~$\alpha$ an
  \emph{$(ab)$-standard colouring} if $\alpha(x_2)=a$ and $\alpha(x_5)=b$.
  Easy counting shows that for fixed $a,b$, there are~15 $(ab)$-standard
  colourings. (\,There are~2 choices for each of the other~4 vertices, but
  one of the resulting~16 3-colourings is not strong.\,) As there are~6
  choices for pairs $a,b$ from~3 colours, there are a total of~90 strong
  3-colourings of~$I$.

  We will prove the proposition by combining the following two steps.

  \medskip\noindent
  \textbf{Step 1}:\quad For given~$a,b$, there is a path containing all
  $(ab)$-standard colourings.

  \begin{center}
    \unitlength0.20mm
    \begin{picture}(80,120)(5,-30)
      \put(25,0){\circle*{5}}\put(25,40){\circle*{5}}
      \put(45,0){\circle*{5}}\put(45,40){\circle*{5}}
      \put(65,0){\circle*{5}}\put(65,40){\circle*{5}}
      \put(25,0){\line(1,0){40}}
      \put(25,40){\line(1,0){40}}
      \put(45,0){\line(0,1){40}}
      \put(25,55){\makebox(0,0)[b]{$b$}}
      \put(45,55){\makebox(0,0)[b]{$a$}}
      \put(65,55){\makebox(0,0)[b]{$b$}}
      \put(25,-25){\makebox(0,0)[b]{$c$}}
      \put(45,-25){\makebox(0,0)[b]{$b$}}
      \put(65,-25){\makebox(0,0)[b]{$a$}}
    \end{picture}
    \begin{picture}(25,120)(0,-30)
      \put(0,20){\vector(1,0){25}}
    \end{picture}
    \begin{picture}(80,120)(5,-30)
      \put(25,0){\circle*{5}}\put(25,40){\circle*{5}}
      \put(45,0){\circle*{5}}\put(45,40){\circle*{5}}
      \put(65,0){\circle*{5}}\put(65,40){\circle*{5}}
      \put(25,0){\line(1,0){40}}
      \put(25,40){\line(1,0){40}}
      \put(45,0){\line(0,1){40}}
      \put(25,55){\makebox(0,0)[b]{$b$}}
      \put(45,55){\makebox(0,0)[b]{$a$}}
      \put(65,55){\makebox(0,0)[b]{$c$}}
      \put(25,-25){\makebox(0,0)[b]{$c$}}
      \put(45,-25){\makebox(0,0)[b]{$b$}}
      \put(65,-25){\makebox(0,0)[b]{$a$}}
    \end{picture}
    \begin{picture}(25,120)(0,-30)
      \put(0,20){\vector(1,0){25}}
    \end{picture}
    \begin{picture}(80,120)(5,-30)
      \put(25,0){\circle*{5}}\put(25,40){\circle*{5}}
      \put(45,0){\circle*{5}}\put(45,40){\circle*{5}}
      \put(65,0){\circle*{5}}\put(65,40){\circle*{5}}
      \put(25,0){\line(1,0){40}}
      \put(25,40){\line(1,0){40}}
      \put(45,0){\line(0,1){40}}
      \put(25,55){\makebox(0,0)[b]{$b$}}
      \put(45,55){\makebox(0,0)[b]{$a$}}
      \put(65,55){\makebox(0,0)[b]{$c$}}
      \put(25,-25){\makebox(0,0)[b]{$a$}}
      \put(45,-25){\makebox(0,0)[b]{$b$}}
      \put(65,-25){\makebox(0,0)[b]{$a$}}
    \end{picture}
    \begin{picture}(25,120)(0,-30)
      \put(0,20){\vector(1,0){25}}
    \end{picture}
    \begin{picture}(80,120)(5,-30)
      \put(25,0){\circle*{5}}\put(25,40){\circle*{5}}
      \put(45,0){\circle*{5}}\put(45,40){\circle*{5}}
      \put(65,0){\circle*{5}}\put(65,40){\circle*{5}}
      \put(25,0){\line(1,0){40}}
      \put(25,40){\line(1,0){40}}
      \put(45,0){\line(0,1){40}}
      \put(25,55){\makebox(0,0)[b]{$c$}}
      \put(45,55){\makebox(0,0)[b]{$a$}}
      \put(65,55){\makebox(0,0)[b]{$c$}}
      \put(25,-25){\makebox(0,0)[b]{$a$}}
      \put(45,-25){\makebox(0,0)[b]{$b$}}
      \put(65,-25){\makebox(0,0)[b]{$a$}}
    \end{picture}
    \begin{picture}(25,120)(0,-30)
      \put(0,20){\vector(1,0){25}}
    \end{picture}
    \begin{picture}(80,120)(5,-30)
      \put(25,0){\circle*{5}}\put(25,40){\circle*{5}}
      \put(45,0){\circle*{5}}\put(45,40){\circle*{5}}
      \put(65,0){\circle*{5}}\put(65,40){\circle*{5}}
      \put(25,0){\line(1,0){40}}
      \put(25,40){\line(1,0){40}}
      \put(45,0){\line(0,1){40}}
      \put(25,55){\makebox(0,0)[b]{$c$}}
      \put(45,55){\makebox(0,0)[b]{$a$}}
      \put(65,55){\makebox(0,0)[b]{$c$}}
      \put(25,-25){\makebox(0,0)[b]{$a$}}
      \put(45,-25){\makebox(0,0)[b]{$b$}}
      \put(65,-25){\makebox(0,0)[b]{$c$}}
    \end{picture}
    \begin{picture}(25,120)(0,-30)
      \put(0,20){\vector(1,0){25}}
    \end{picture}
    \begin{picture}(80,120)(5,-30)
      \put(25,0){\circle*{5}}\put(25,40){\circle*{5}}
      \put(45,0){\circle*{5}}\put(45,40){\circle*{5}}
      \put(65,0){\circle*{5}}\put(65,40){\circle*{5}}
      \put(25,0){\line(1,0){40}}
      \put(25,40){\line(1,0){40}}
      \put(45,0){\line(0,1){40}}
      \put(25,55){\makebox(0,0)[b]{$b$}}
      \put(45,55){\makebox(0,0)[b]{$a$}}
      \put(65,55){\makebox(0,0)[b]{$c$}}
      \put(25,-25){\makebox(0,0)[b]{$a$}}
      \put(45,-25){\makebox(0,0)[b]{$b$}}
      \put(65,-25){\makebox(0,0)[b]{$c$}}
    \end{picture}
    \begin{picture}(25,120)(0,-30)
      \put(0,20){\vector(1,0){25}}
    \end{picture}
    \begin{picture}(80,120)(5,-30)
      \put(25,0){\circle*{5}}\put(25,40){\circle*{5}}
      \put(45,0){\circle*{5}}\put(45,40){\circle*{5}}
      \put(65,0){\circle*{5}}\put(65,40){\circle*{5}}
      \put(25,0){\line(1,0){40}}
      \put(25,40){\line(1,0){40}}
      \put(45,0){\line(0,1){40}}
      \put(25,55){\makebox(0,0)[b]{$b$}}
      \put(45,55){\makebox(0,0)[b]{$a$}}
      \put(65,55){\makebox(0,0)[b]{$c$}}
      \put(25,-25){\makebox(0,0)[b]{$c$}}
      \put(45,-25){\makebox(0,0)[b]{$b$}}
      \put(65,-25){\makebox(0,0)[b]{$c$}}
    \end{picture}
    \begin{picture}(25,120)(0,-30)
      \put(0,20){\vector(1,0){25}}
    \end{picture}
    \begin{picture}(80,120)(5,-30)
      \put(25,0){\circle*{5}}\put(25,40){\circle*{5}}
      \put(45,0){\circle*{5}}\put(45,40){\circle*{5}}
      \put(65,0){\circle*{5}}\put(65,40){\circle*{5}}
      \put(25,0){\line(1,0){40}}
      \put(25,40){\line(1,0){40}}
      \put(45,0){\line(0,1){40}}
      \put(25,55){\makebox(0,0)[b]{$c$}}
      \put(45,55){\makebox(0,0)[b]{$a$}}
      \put(65,55){\makebox(0,0)[b]{$c$}}
      \put(25,-25){\makebox(0,0)[b]{$c$}}
      \put(45,-25){\makebox(0,0)[b]{$b$}}
      \put(65,-25){\makebox(0,0)[b]{$c$}}
    \end{picture}
    \begin{picture}(25,120)(0,-30)
      \put(0,20){\vector(1,0){25}}
    \end{picture}
    \begin{picture}(80,120)(5,-30)
      \put(25,0){\circle*{5}}\put(25,40){\circle*{5}}
      \put(45,0){\circle*{5}}\put(45,40){\circle*{5}}
      \put(65,0){\circle*{5}}\put(65,40){\circle*{5}}
      \put(25,0){\line(1,0){40}}
      \put(25,40){\line(1,0){40}}
      \put(45,0){\line(0,1){40}}
      \put(25,55){\makebox(0,0)[b]{$c$}}
      \put(45,55){\makebox(0,0)[b]{$a$}}
      \put(65,55){\makebox(0,0)[b]{$c$}}
      \put(25,-25){\makebox(0,0)[b]{$c$}}
      \put(45,-25){\makebox(0,0)[b]{$b$}}
      \put(65,-25){\makebox(0,0)[b]{$a$}}
    \end{picture}
    \begin{picture}(25,120)(0,-30)
      \put(0,20){\vector(1,0){25}}
    \end{picture}
    \begin{picture}(80,120)(5,-30)
      \put(25,0){\circle*{5}}\put(25,40){\circle*{5}}
      \put(45,0){\circle*{5}}\put(45,40){\circle*{5}}
      \put(65,0){\circle*{5}}\put(65,40){\circle*{5}}
      \put(25,0){\line(1,0){40}}
      \put(25,40){\line(1,0){40}}
      \put(45,0){\line(0,1){40}}
      \put(25,55){\makebox(0,0)[b]{$c$}}
      \put(45,55){\makebox(0,0)[b]{$a$}}
      \put(65,55){\makebox(0,0)[b]{$b$}}
      \put(25,-25){\makebox(0,0)[b]{$c$}}
      \put(45,-25){\makebox(0,0)[b]{$b$}}
      \put(65,-25){\makebox(0,0)[b]{$a$}}
    \end{picture}
    \begin{picture}(25,120)(0,-30)
      \put(0,20){\vector(1,0){25}}
    \end{picture}
    \begin{picture}(80,120)(5,-30)
      \put(25,0){\circle*{5}}\put(25,40){\circle*{5}}
      \put(45,0){\circle*{5}}\put(45,40){\circle*{5}}
      \put(65,0){\circle*{5}}\put(65,40){\circle*{5}}
      \put(25,0){\line(1,0){40}}
      \put(25,40){\line(1,0){40}}
      \put(45,0){\line(0,1){40}}
      \put(25,55){\makebox(0,0)[b]{$c$}}
      \put(45,55){\makebox(0,0)[b]{$a$}}
      \put(65,55){\makebox(0,0)[b]{$b$}}
      \put(25,-25){\makebox(0,0)[b]{$c$}}
      \put(45,-25){\makebox(0,0)[b]{$b$}}
      \put(65,-25){\makebox(0,0)[b]{$c$}}
    \end{picture}
    \begin{picture}(25,120)(0,-30)
      \put(0,20){\vector(1,0){25}}
    \end{picture}
    \begin{picture}(80,120)(5,-30)
      \put(25,0){\circle*{5}}\put(25,40){\circle*{5}}
      \put(45,0){\circle*{5}}\put(45,40){\circle*{5}}
      \put(65,0){\circle*{5}}\put(65,40){\circle*{5}}
      \put(25,0){\line(1,0){40}}
      \put(25,40){\line(1,0){40}}
      \put(45,0){\line(0,1){40}}
      \put(25,55){\makebox(0,0)[b]{$b$}}
      \put(45,55){\makebox(0,0)[b]{$a$}}
      \put(65,55){\makebox(0,0)[b]{$b$}}
      \put(25,-25){\makebox(0,0)[b]{$c$}}
      \put(45,-25){\makebox(0,0)[b]{$b$}}
      \put(65,-25){\makebox(0,0)[b]{$c$}}
    \end{picture}
    \begin{picture}(25,120)(0,-30)
      \put(0,20){\vector(1,0){25}}
    \end{picture}
    \begin{picture}(80,120)(5,-30)
      \put(25,0){\circle*{5}}\put(25,40){\circle*{5}}
      \put(45,0){\circle*{5}}\put(45,40){\circle*{5}}
      \put(65,0){\circle*{5}}\put(65,40){\circle*{5}}
      \put(25,0){\line(1,0){40}}
      \put(25,40){\line(1,0){40}}
      \put(45,0){\line(0,1){40}}
      \put(25,55){\makebox(0,0)[b]{$b$}}
      \put(45,55){\makebox(0,0)[b]{$a$}}
      \put(65,55){\makebox(0,0)[b]{$b$}}
      \put(25,-25){\makebox(0,0)[b]{$a$}}
      \put(45,-25){\makebox(0,0)[b]{$b$}}
      \put(65,-25){\makebox(0,0)[b]{$c$}}
    \end{picture}
    \begin{picture}(25,120)(0,-30)
      \put(0,20){\vector(1,0){25}}
    \end{picture}
    \begin{picture}(80,120)(5,-30)
      \put(25,0){\circle*{5}}\put(25,40){\circle*{5}}
      \put(45,0){\circle*{5}}\put(45,40){\circle*{5}}
      \put(65,0){\circle*{5}}\put(65,40){\circle*{5}}
      \put(25,0){\line(1,0){40}}
      \put(25,40){\line(1,0){40}}
      \put(45,0){\line(0,1){40}}
      \put(25,55){\makebox(0,0)[b]{$c$}}
      \put(45,55){\makebox(0,0)[b]{$a$}}
      \put(65,55){\makebox(0,0)[b]{$b$}}
      \put(25,-25){\makebox(0,0)[b]{$a$}}
      \put(45,-25){\makebox(0,0)[b]{$b$}}
      \put(65,-25){\makebox(0,0)[b]{$c$}}
    \end{picture}
    \begin{picture}(25,120)(0,-30)
      \put(0,20){\vector(1,0){25}}
    \end{picture}
    \begin{picture}(80,120)(5,-30)
      \put(25,0){\circle*{5}}\put(25,40){\circle*{5}}
      \put(45,0){\circle*{5}}\put(45,40){\circle*{5}}
      \put(65,0){\circle*{5}}\put(65,40){\circle*{5}}
      \put(25,0){\line(1,0){40}}
      \put(25,40){\line(1,0){40}}
      \put(45,0){\line(0,1){40}}
      \put(25,55){\makebox(0,0)[b]{$c$}}
      \put(45,55){\makebox(0,0)[b]{$a$}}
      \put(65,55){\makebox(0,0)[b]{$b$}}
      \put(25,-25){\makebox(0,0)[b]{$a$}}
      \put(45,-25){\makebox(0,0)[b]{$b$}}
      \put(65,-25){\makebox(0,0)[b]{$a$}}
    \end{picture}
  \end{center}

  \medskip\noindent
  \textbf{Step 2}:\quad There is a path containing at least one colouring
  from each type of standard colourings.

  \begin{center}
    \unitlength0.20mm
    \begin{picture}(80,120)(5,-30)
      \put(25,0){\circle*{5}}\put(25,40){\circle*{5}}
      \put(45,0){\circle*{5}}\put(45,40){\circle*{5}}
      \put(65,0){\circle*{5}}\put(65,40){\circle*{5}}
      \put(25,0){\line(1,0){40}}
      \put(25,40){\line(1,0){40}}
      \put(45,0){\line(0,1){40}}
      \put(25,55){\makebox(0,0)[b]{$c$}}
      \put(45,55){\makebox(0,0)[b]{$b$}}
      \put(65,55){\makebox(0,0)[b]{$a$}}
      \put(25,-25){\makebox(0,0)[b]{$b$}}
      \put(45,-25){\makebox(0,0)[b]{$c$}}
      \put(65,-25){\makebox(0,0)[b]{$b$}}
    \end{picture}
    \begin{picture}(25,120)(0,-30)
      \put(0,20){\vector(1,0){25}}
    \end{picture}
    \begin{picture}(80,120)(5,-30)
      \put(25,0){\circle*{5}}\put(25,40){\circle*{5}}
      \put(45,0){\circle*{5}}\put(45,40){\circle*{5}}
      \put(65,0){\circle*{5}}\put(65,40){\circle*{5}}
      \put(25,0){\line(1,0){40}}
      \put(25,40){\line(1,0){40}}
      \put(45,0){\line(0,1){40}}
      \put(25,55){\makebox(0,0)[b]{$c$}}
      \put(45,55){\makebox(0,0)[b]{$b$}}
      \put(65,55){\makebox(0,0)[b]{$a$}}
      \put(25,-25){\makebox(0,0)[b]{$b$}}
      \put(45,-25){\makebox(0,0)[b]{$a$}}
      \put(65,-25){\makebox(0,0)[b]{$b$}}
    \end{picture}
    \begin{picture}(25,120)(0,-30)
      \put(0,20){\vector(1,0){25}}
    \end{picture}
    \begin{picture}(80,120)(5,-30)
      \put(25,0){\circle*{5}}\put(25,40){\circle*{5}}
      \put(45,0){\circle*{5}}\put(45,40){\circle*{5}}
      \put(65,0){\circle*{5}}\put(65,40){\circle*{5}}
      \put(25,0){\line(1,0){40}}
      \put(25,40){\line(1,0){40}}
      \put(45,0){\line(0,1){40}}
      \put(25,55){\makebox(0,0)[b]{$c$}}
      \put(45,55){\makebox(0,0)[b]{$b$}}
      \put(65,55){\makebox(0,0)[b]{$a$}}
      \put(25,-25){\makebox(0,0)[b]{$c$}}
      \put(45,-25){\makebox(0,0)[b]{$a$}}
      \put(65,-25){\makebox(0,0)[b]{$b$}}
    \end{picture}
    \begin{picture}(25,120)(0,-30)
      \put(0,20){\vector(1,0){25}}
    \end{picture}
    \begin{picture}(80,120)(5,-30)
      \put(25,0){\circle*{5}}\put(25,40){\circle*{5}}
      \put(45,0){\circle*{5}}\put(45,40){\circle*{5}}
      \put(65,0){\circle*{5}}\put(65,40){\circle*{5}}
      \put(25,0){\line(1,0){40}}
      \put(25,40){\line(1,0){40}}
      \put(45,0){\line(0,1){40}}
      \put(25,55){\makebox(0,0)[b]{$a$}}
      \put(45,55){\makebox(0,0)[b]{$b$}}
      \put(65,55){\makebox(0,0)[b]{$a$}}
      \put(25,-25){\makebox(0,0)[b]{$c$}}
      \put(45,-25){\makebox(0,0)[b]{$a$}}
      \put(65,-25){\makebox(0,0)[b]{$b$}}
    \end{picture}
    \begin{picture}(25,120)(0,-30)
      \put(0,20){\vector(1,0){25}}
    \end{picture}
    \begin{picture}(80,120)(5,-30)
      \put(25,0){\circle*{5}}\put(25,40){\circle*{5}}
      \put(45,0){\circle*{5}}\put(45,40){\circle*{5}}
      \put(65,0){\circle*{5}}\put(65,40){\circle*{5}}
      \put(25,0){\line(1,0){40}}
      \put(25,40){\line(1,0){40}}
      \put(45,0){\line(0,1){40}}
      \put(25,55){\makebox(0,0)[b]{$a$}}
      \put(45,55){\makebox(0,0)[b]{$c$}}
      \put(65,55){\makebox(0,0)[b]{$a$}}
      \put(25,-25){\makebox(0,0)[b]{$c$}}
      \put(45,-25){\makebox(0,0)[b]{$a$}}
      \put(65,-25){\makebox(0,0)[b]{$b$}}
    \end{picture}
    \begin{picture}(25,120)(0,-30)
      \put(0,20){\vector(1,0){25}}
    \end{picture}
    \begin{picture}(80,120)(5,-30)
      \put(25,0){\circle*{5}}\put(25,40){\circle*{5}}
      \put(45,0){\circle*{5}}\put(45,40){\circle*{5}}
      \put(65,0){\circle*{5}}\put(65,40){\circle*{5}}
      \put(25,0){\line(1,0){40}}
      \put(25,40){\line(1,0){40}}
      \put(45,0){\line(0,1){40}}
      \put(25,55){\makebox(0,0)[b]{$a$}}
      \put(45,55){\makebox(0,0)[b]{$c$}}
      \put(65,55){\makebox(0,0)[b]{$b$}}
      \put(25,-25){\makebox(0,0)[b]{$c$}}
      \put(45,-25){\makebox(0,0)[b]{$a$}}
      \put(65,-25){\makebox(0,0)[b]{$b$}}
    \end{picture}
    \begin{picture}(25,120)(0,-30)
      \put(0,20){\vector(1,0){25}}
    \end{picture}
    \begin{picture}(80,120)(5,-30)
      \put(25,0){\circle*{5}}\put(25,40){\circle*{5}}
      \put(45,0){\circle*{5}}\put(45,40){\circle*{5}}
      \put(65,0){\circle*{5}}\put(65,40){\circle*{5}}
      \put(25,0){\line(1,0){40}}
      \put(25,40){\line(1,0){40}}
      \put(45,0){\line(0,1){40}}
      \put(25,55){\makebox(0,0)[b]{$a$}}
      \put(45,55){\makebox(0,0)[b]{$c$}}
      \put(65,55){\makebox(0,0)[b]{$b$}}
      \put(25,-25){\makebox(0,0)[b]{$c$}}
      \put(45,-25){\makebox(0,0)[b]{$a$}}
      \put(65,-25){\makebox(0,0)[b]{$c$}}
    \end{picture}
    \begin{picture}(25,120)(0,-30)
      \put(0,20){\vector(1,0){25}}
    \end{picture}
    \begin{picture}(80,120)(5,-30)
      \put(25,0){\circle*{5}}\put(25,40){\circle*{5}}
      \put(45,0){\circle*{5}}\put(45,40){\circle*{5}}
      \put(65,0){\circle*{5}}\put(65,40){\circle*{5}}
      \put(25,0){\line(1,0){40}}
      \put(25,40){\line(1,0){40}}
      \put(45,0){\line(0,1){40}}
      \put(25,55){\makebox(0,0)[b]{$a$}}
      \put(45,55){\makebox(0,0)[b]{$c$}}
      \put(65,55){\makebox(0,0)[b]{$b$}}
      \put(25,-25){\makebox(0,0)[b]{$c$}}
      \put(45,-25){\makebox(0,0)[b]{$b$}}
      \put(65,-25){\makebox(0,0)[b]{$c$}}
    \end{picture}
    \begin{picture}(25,120)(0,-30)
      \put(0,20){\vector(1,0){25}}
    \end{picture}
    \begin{picture}(80,120)(5,-30)
      \put(25,0){\circle*{5}}\put(25,40){\circle*{5}}
      \put(45,0){\circle*{5}}\put(45,40){\circle*{5}}
      \put(65,0){\circle*{5}}\put(65,40){\circle*{5}}
      \put(25,0){\line(1,0){40}}
      \put(25,40){\line(1,0){40}}
      \put(45,0){\line(0,1){40}}
      \put(25,55){\makebox(0,0)[b]{$a$}}
      \put(45,55){\makebox(0,0)[b]{$c$}}
      \put(65,55){\makebox(0,0)[b]{$b$}}
      \put(25,-25){\makebox(0,0)[b]{$a$}}
      \put(45,-25){\makebox(0,0)[b]{$b$}}
      \put(65,-25){\makebox(0,0)[b]{$c$}}
    \end{picture}
    \begin{picture}(25,120)(0,-30)
      \put(0,20){\vector(1,0){25}}
    \end{picture}
    \begin{picture}(80,120)(5,-30)
      \put(25,0){\circle*{5}}\put(25,40){\circle*{5}}
      \put(45,0){\circle*{5}}\put(45,40){\circle*{5}}
      \put(65,0){\circle*{5}}\put(65,40){\circle*{5}}
      \put(25,0){\line(1,0){40}}
      \put(25,40){\line(1,0){40}}
      \put(45,0){\line(0,1){40}}
      \put(25,55){\makebox(0,0)[b]{$b$}}
      \put(45,55){\makebox(0,0)[b]{$c$}}
      \put(65,55){\makebox(0,0)[b]{$b$}}
      \put(25,-25){\makebox(0,0)[b]{$a$}}
      \put(45,-25){\makebox(0,0)[b]{$b$}}
      \put(65,-25){\makebox(0,0)[b]{$c$}}
    \end{picture}
    \begin{picture}(25,120)(0,-30)
      \put(0,20){\vector(1,0){25}}
    \end{picture}
    \begin{picture}(80,120)(5,-30)
      \put(25,0){\circle*{5}}\put(25,40){\circle*{5}}
      \put(45,0){\circle*{5}}\put(45,40){\circle*{5}}
      \put(65,0){\circle*{5}}\put(65,40){\circle*{5}}
      \put(25,0){\line(1,0){40}}
      \put(25,40){\line(1,0){40}}
      \put(45,0){\line(0,1){40}}
      \put(25,55){\makebox(0,0)[b]{$b$}}
      \put(45,55){\makebox(0,0)[b]{$a$}}
      \put(65,55){\makebox(0,0)[b]{$b$}}
      \put(25,-25){\makebox(0,0)[b]{$a$}}
      \put(45,-25){\makebox(0,0)[b]{$b$}}
      \put(65,-25){\makebox(0,0)[b]{$c$}}
    \end{picture}
    \begin{picture}(25,120)(0,-30)
      \put(0,20){\vector(1,0){25}}
    \end{picture}
    \begin{picture}(80,120)(5,-30)
      \put(25,0){\circle*{5}}\put(25,40){\circle*{5}}
      \put(45,0){\circle*{5}}\put(45,40){\circle*{5}}
      \put(65,0){\circle*{5}}\put(65,40){\circle*{5}}
      \put(25,0){\line(1,0){40}}
      \put(25,40){\line(1,0){40}}
      \put(45,0){\line(0,1){40}}
      \put(25,55){\makebox(0,0)[b]{$b$}}
      \put(45,55){\makebox(0,0)[b]{$a$}}
      \put(65,55){\makebox(0,0)[b]{$c$}}
      \put(25,-25){\makebox(0,0)[b]{$a$}}
      \put(45,-25){\makebox(0,0)[b]{$b$}}
      \put(65,-25){\makebox(0,0)[b]{$c$}}
    \end{picture}
    \begin{picture}(25,120)(0,-30)
      \put(0,20){\vector(1,0){25}}
    \end{picture}
    \begin{picture}(80,120)(5,-30)
      \put(25,0){\circle*{5}}\put(25,40){\circle*{5}}
      \put(45,0){\circle*{5}}\put(45,40){\circle*{5}}
      \put(65,0){\circle*{5}}\put(65,40){\circle*{5}}
      \put(25,0){\line(1,0){40}}
      \put(25,40){\line(1,0){40}}
      \put(45,0){\line(0,1){40}}
      \put(25,55){\makebox(0,0)[b]{$b$}}
      \put(45,55){\makebox(0,0)[b]{$a$}}
      \put(65,55){\makebox(0,0)[b]{$c$}}
      \put(25,-25){\makebox(0,0)[b]{$a$}}
      \put(45,-25){\makebox(0,0)[b]{$b$}}
      \put(65,-25){\makebox(0,0)[b]{$a$}}
    \end{picture}
    \begin{picture}(25,120)(0,-30)
      \put(0,20){\vector(1,0){25}}
    \end{picture}
    \begin{picture}(80,120)(5,-30)
      \put(25,0){\circle*{5}}\put(25,40){\circle*{5}}
      \put(45,0){\circle*{5}}\put(45,40){\circle*{5}}
      \put(65,0){\circle*{5}}\put(65,40){\circle*{5}}
      \put(25,0){\line(1,0){40}}
      \put(25,40){\line(1,0){40}}
      \put(45,0){\line(0,1){40}}
      \put(25,55){\makebox(0,0)[b]{$b$}}
      \put(45,55){\makebox(0,0)[b]{$a$}}
      \put(65,55){\makebox(0,0)[b]{$c$}}
      \put(25,-25){\makebox(0,0)[b]{$a$}}
      \put(45,-25){\makebox(0,0)[b]{$c$}}
      \put(65,-25){\makebox(0,0)[b]{$a$}}
    \end{picture}
  \end{center}

  \medskip\noindent
  These two steps, together with appropriate renaming of the colours, will
  give all that is needed to transform any strong 3-colouring of~$I$ into
  any other strong 3-colouring.
\end{proof}

\begin{thm}\label{D.10}
  Let~$T$ be a tree. Then $S_3(T)$ is connected if and only if~$T$
  contains~$P_5$ or~$I$ as a subgraph.
\end{thm}

\begin{proof}
  Since $S_3(P_5)$ and $S_3(I)$ are connected, one direction is immediately
  proved by using Theorem~\ref{Z.6}.

  For the other direction, suppose that~$T$ does not contain~$P_5$,
  nor~$I$. Since~$P_5$ is a path with~4 edges, the longest path in~$T$ can
  have length at most~3. Thus $T$ has to be one of the following\,: $K_1$,
  $P_2$, $K_{1,m}$, $m\ge2$, or~$\Psi_n$, $n\ge1$. Note that~$T$ cannot be
  $\Phi_{p,q}$ for all $p,q\ge2$ since then it would contain~$I$ as a
  subgraph. Since~$K_1$ and~$P_2$ have fewer than~3 vertices, $T$ cannot be
  one of these two graphs. We already saw that $S_3(K_{1,m})$, $m\ge2$, and
  $S_3(\Psi_n)$, $n\ge1$, are disconnected. We can conclude that $S_3(T)$
  is not connected, which completes the proof.
\end{proof}

\section{Discussion}\label{sec6}

  We realise that this note contains only some first results on the strong
  colour graphs. In comparison, there is a growing body of literature on the
  connectivity of the normal colour graph\,:
  \cite{vertex-colouring,mixing,YMCAL,OMPH}. An interesting direction of
  future research would be to investigate how far the theory of strong colour
  graphs can be reduced to the theory of normal colour graphs.

  We have already seen in Lemma~\ref{Z.10} that if the strong colour graph
  $S_k(G)$ is connected for some~$G$ and~$k$, then so is the normal colour
  graph~$C_k(G)$. In general, the reverse direction is not true. For
  instance, for all $m,n\ge2$, for the complete bipartite graphs~$K_{m,n}$ we
  have that for $k\ge3$, $C_k(K_{m,n})$ is connected, whereas $S_k(K_{m,n})$
  is not connected. In fact, we've already seen that if~$G$ has a complete
  bipartite graph as a spanning subgraph, then $S_k(G)$ is never connected.
  For $k\ge3$ it is not hard to construct other graphs apart from complete
  bipartite graphs that have this property, but all examples we know of have
  a fairly special structure (\,see for instance the trees~$\Psi_n$ in
  Figure~\ref{f5}\,). This makes us raise the following question.

\begin{ques}\label{Q1}
  Is it possible to completely describe a class of graphs~$\mathcal{H}$
  so that if~$G$ does not contain a graph from~$\mathcal{H}$ as a spanning
  subgraph, then~$C_k(G)$ is connected if and only if~$S_k(G)$ is
  connected\,?
\end{ques}

\end{document}